\RequirePackage{fix-cm}
\documentclass{article}                     
\usepackage[margin=1in]{geometry}
\usepackage[ansinew]{inputenc}
\usepackage[english]{babel}
\usepackage{fancyhdr}
\usepackage[active]{srcltx}
\usepackage{amsmath,amssymb}
\usepackage{xcolor}
\usepackage{array}
\usepackage{textcomp}
\usepackage{tikz}
\usepackage[ruled,vlined]{algorithm2e}
\usepackage{cuted}
\usepackage{placeins}
\usepackage[normalem]{ulem}
\usepackage{hyperref}
\usepackage{xcolor}
\usepackage{bm}
\usepackage{stackengine}
\usepackage{authblk}

\newcommand{\bsn}[1]{\boldsymbol{#1}}
\newcommand{\fun}[1]{\bm{\mathit{#1}}}

\newcommand{\p}[2][p]{{#2^{(#1)}}}
\renewcommand{\P}[2][p]{{\left[#2\right]_{#1}}}

\def\delequal{\mathrel{\ensurestackMath{\stackon[1pt]{=}{\scriptscriptstyle \Delta}}}}

\definecolor{cpurple}{RGB}{163,99,217}
\definecolor{cred}{RGB}{219,56,56}
\definecolor{corange}{RGB}{249,162,40} 
\definecolor{cgreen}{RGB}{178,194,37} 
\definecolor{cblue}{RGB}{64,164,216} 

\newcommand{\sqboxs}{1.4ex}
\newcommand{\sqboxpurple}[1]{\textcolor{cpurple}{\rule{\sqboxs}{\sqboxs}}}
\newcommand{\sqboxred}[1]{\textcolor{cred}{\rule{\sqboxs}{\sqboxs}}}
\newcommand{\sqboxorange}[1]{\textcolor{corange}{\rule{\sqboxs}{\sqboxs}}}
\newcommand{\sqboxgreen}[1]{\textcolor{cgreen}{\rule{\sqboxs}{\sqboxs}}}
\newcommand{\sqboxblue}[1]{\textcolor{cblue}{\rule{\sqboxs}{\sqboxs}}}

\newcommand{\veps}{\varepsilon}

\newcommand{\hFu}{\hat{\mathbf{F}}}
\newcommand{\hn}{\hat{n}}

\newcommand{\hlm}{\hat{\lambda}}

\newcommand{\tFu}{\tilde{\mathbf{F}}}
\newcommand{\tKu}{\tilde{\mathbf{K}}}

\newcommand{\sym}{\mathrm{sym}}
\newcommand{\dd}{\mathrm{d}}

\newcommand{\Bu}{\bsn{B}}

\newcommand{\Eu}{\bsn{E}}

\newcommand{\Pu}{\bsn{P}}

\newcommand{\Su}{\bsn{S}}

\newcommand{\Wu}{\mathbf{W}}

\newcommand{\eu}{\bsn{e}}

\newcommand{\pu}{\bsn{p}}

\newcommand{\lu}{\bsn{l}}
\newcommand{\nuu}{\bsn{n}}
\newcommand{\uu}{\bsn{u}}

\newcommand{\wu}{\bsn{w}}

\newcommand{\cA}{\mathcal{A}}
\newcommand{\cC}{\mathcal{C}}
\newcommand{\cH}{\mathcal{H}}
\newcommand{\cO}{\mathcal{O}}

\newcommand{\cQ}{\mathcal{Q}}

\newcommand{\Omegau}{\bsn{\Omega}}

\newcommand{\Uu}{\mathbf{U}}

\newcommand{\tUu}{\tilde{\mathbf{U}}}
\newcommand{\tVu}{\tilde{\mathbf{V}}}

\newcommand{\Mu}{\mathbf{M}}
\newcommand{\Fu}{\mathbf{F}}
\newcommand{\Ku}{\mathbf{K}}
\newcommand{\Cu}{\mathbf{C}}
\newcommand{\G}{\fun{G}}
\newcommand{\tG}{\tilde{\fun{G}}}
\renewcommand{\H}{\fun{H}}

\newcommand{\bX}{\mathbb{X}}
\newcommand{\bY}{\mathbb{Y}}
\newcommand{\bPhi}{{\mathbf{\Phi}}}

\newcommand{\Xu}{\mathbf{X}}
\newcommand{\Yu}{\mathbf{Y}}
\newcommand{\phiu}{\pmb{\phi}}
\newcommand{\lambdau}{\bsn{\lambda}}

\newcommand{\FU}{\bm{\mu}}
\newcommand{\FV}{\bm{\nu}}
\newcommand{\hFU}{\hat{\bm{\mu}}}
\newcommand{\hFV}{\hat{\bm{\nu}}}
\newcommand{\gFU}{\mathring{\bm{\mu}}}

\newcommand{\FG}{\mathbf{G}}
\newcommand{\FH}{\mathbf{H}}
\newcommand{\hFG}{\hat{\mathbf{G}}}
\newcommand{\hFH}{\hat{\mathbf{H}}}

\newcommand{\hsigma}{\hat{\sigma}}

\newcommand{\gsigma}{\mathring{\sigma}}

\newcommand{\gE}{\mathring{\mathbf{E}}}

\newcommand{\WU}{\mathbf{\Psi}}
\newcommand{\WV}{\mathbf{\Upsilon}}
\newcommand{\hWU}{\hat{\mathbf{\Psi}}}
\newcommand{\hWV}{\hat{\mathbf{\Upsilon}}}
\newcommand{\gWU}{\mathring{\mathbf{\Psi}}}
\newcommand{\gWV}{\mathring{\mathbf{\Upsilon}}}

\newcommand{\WUfun}{\fun{\Psi}(\zu)}
\newcommand{\WVfun}{\fun{\Upsilon}(\zu)}
\newcommand{\hWUfun}{\hat{\fun{\Psi}}(\zu)}
\newcommand{\hWVfun}{\hat{\fun{\Upsilon}}(\zu)}

\newcommand{\ffun}{\fun{f}(\zu)}
\newcommand{\hffun}{\hat{\fun{f}}(\zu)}

\newcommand{\WUfunsim}{\fun{\Psi}}
\newcommand{\WVfunsim}{\fun{\Upsilon}}
\newcommand{\hWUfunsim}{\hat{\fun{{\Psi}}}}
\newcommand{\hWVfunsim}{\hat{\fun{{\Upsilon}}}}

\newcommand{\ffunsim}{\fun{{f}}}
\newcommand{\hffunsim}{\hat{\fun{{f}}}}

\newcommand{\zerou}{\mathbf{0}}

\newcommand{\zr}{z}
\newcommand{\fr}{f}
\newcommand{\zu}{\mathbf{z}}
\newcommand{\fu}{\mathbf{f}}

\newcommand{\hf}{\hat{{f}}}

\newcommand{\gf}{\mathring{{f}}}

\newcommand{\Is}{\mathcal{I}}

\DeclareMathAlphabet{\mymathbb}{U}{BOONDOX-ds}{m}{n}
\newcommand{\iu}{\mymathbb{i}}

\definecolor{myred}{rgb}{0.9,0,0}
\definecolor{Color_FS}{rgb}{.7,0,.7}
\definecolor{Color_NF}{rgb}{0,0,.7}
\definecolor{Color_NFTO}{rgb}{0,.5,.2}
\definecolor{Color_MD}{rgb}{.75,.35,0}
\definecolor{Color_SMD}{rgb}{1,.7,0}

\begin{document}

\title{Nonlinear model order reduction of resonant piezoelectric micro-actuators: an invariant manifold approach}

\author[1]{Andrea Opreni}
\author[1]{Giorgio Gobat\footnote{giorgio.gobat@polimi.it}}
\author[2]{Cyril Touz\'e}
\author[1]{Attilio Frangi}
\affil[1]{Department of Civil and Environmental Engineering, Politecnico di Milano, P.za Leonardo da Vinci 32, 20133 Milano, Italy}
\affil[2]{Institute of Mechanical Sciences and Industrial Applications (IMSIA)\\ ENSTA Paris - CNRS - EDF - CEA - Institut Polytechnique de Paris 828 boulevard des mar\'echaux 91762 Palaiseau cedex}


\maketitle

\begin{abstract}
This paper presents a novel derivation of the direct parametrisation method for invariant manifolds able to build simulation-free reduced-order models for nonlinear piezoelectric structures, with a particular emphasis on applications to Micro-Electro-Mechanical-Systems.
The constitutive model adopted accounts for the hysteretic and electrostrictive response of the piezoelectric material by resorting to 
the Landau-Devonshire theory of ferroelectrics. 
Results are validated with full-order simulations operated with a harmonic balance finite element method to highlight the reliability of the proposed reduction procedure. Numerical results show a remarkable gain in terms of computing time as a result of the dimensionality reduction process over low dimensional invariant sets. Results are also compared with experimental data to highlight the remarkable benefits of the proposed model order reduction technique.
\end{abstract}

\section{Introduction}
\label{sec:intro}

Piezoelectric Micro Electro Mechanical Systems (MEMS) represent nowadays an important class of devices for both actuation and sensing~\cite{vigna22} 
that are often preferred to their capacitive and magnetic counterparts.
In particular, piezoactuation is a key enabling technology for the next generation
of micromirrors, loudspeakers, piezoelectric ultrasonic transducers~\cite{butt2016,filhol2005,devoe17,massimino2018,rahim2018,yang2018}. 
For these applications, every source of nonlinear behaviour
must be predicted and controlled, since the frequency drift of the resonant mode 
with increasing actuation voltage can be disruptive for the correct functioning of the device.
For instance, structural and piezoelectric material nonlinearities are excited when considering large amplitude vibrations that are routinely reached within the operating range of these devices.
Generally, the growing importance of nonlinear effects in MEMS is stimulating 
intensive research as recently reported in {\em e.g.}~\cite{hajjaj2020linear,shaw21,gobat2021backbone,gobat2021reduced}.
The main aim of this work is thus to consider the derivation of fast and accurate reduced models (ROMs) for nonlinear structures subjected to piezoelectric actuation, with application to MEMS devices. To that purpose, the direct parametrisation method for invariant manifolds will be used and extended in order to properly take into account the new effects provided by the piezoeletric coupling. 

Nonlinear model order reduction techniques have been used for decades for geometrically nonlinear structures, as surveyed for example in~\cite{touze2021review}. They 
rely on a different approach compared to linear projection methods, like {\em e.g.} modal decomposition, Proper Orthogonal Decomposition (POD)~\cite{sampaio2007remarks,amabili2007reduced,gobat2022reduced}, or techniques as the Proper Generalized Decomposition (PGD)~\cite{MEYRAND2019}. 
In this respect, nonlinear normal modes (NNMs) defined as invariant manifolds attached at a fixed point to their linear counterpart, are a powerful tool introduced in the seminal work by Shaw and Pierre~\cite{shaw1991nonlinear}. While the first methods to compute these invariant sets applied either the center manifold theory~\cite{shaw1993normal} or the normal form approach~\cite{touze2004hardening,touze2006nonlinear}, an important step forward has been provided by the parametrisation method of invariant manifolds~\cite{cabre2005overview,haro2016parametrisation}, since the two techniques can be embedded in the same framework.

In recent years two important obstacles hindering the use of nonlinear reduction methods based on invariant manifold theory to large scale systems have been overcome. The first one  is the {\em direct} calculation, allowing one to go  from the physical space (FE degrees-of freedoms, dofs) to a reduced subspace spanned by invariant manifolds. This has been achieved with a normal form approach~\cite{vizzaccaro2021direct,opreni2021model} and with the parametrisation method~\cite{jain2021how}. This step represents a very important achievement since previous methods, {\em e.g.} the ones presented in~\cite{pesheck2002galerkin,touze2006nonlinear,haller2016nnm}, need to express the dynamics in the full modal basis as a starting point, which impedes the applicability of the method to large-scale FE problems. The second important achievement is related to the use of arbitrary order expansions, allowing one to offer automated algorithms with the required accuracy that guarantees convergence. This has been first realized using the parametrisation method from the modal space~\cite{ponsioen2018automated,ponsioen2020model}, and has been pushed forward and implemented through a direct approach in~\cite{jain2021how,vizza21high,opreni22high}. In particular, non-autonomous generic forcing terms can be directly treated in the parametrisation procedure as in~\cite{opreni22high}. Whereas the simple solution of directly adding the modal forcing to the reduced dynamics, as used for instance in~\cite{touze2006nonlinear,vizzaccaro2021direct,vizza21high,jain2021how}, gives accurate results when the loading is colinear with the master mode itself, higher-orders on the forcing appear to be of prime importance in other cases \cite{opreni22high} such as non-modal forcing, forcing orthogonal to the master mode, parametric excitation, and the computation of isolated solutions arising in such cases.

Considering applications to MEMS which address different physics such 
as piezomechanical coupling,
the generation of Reduced Order Models is a complex task in itself
as classical linear reduction approaches \cite{massimino2018,korvink15}
cannot include the effects of geometric nonlinearities.
These have been addressed in 
\cite{lazarus2012finite,givois2020experimental,thomas2013efficient}, where a ROM for coupled piezoelectricity was developed using structural theories 
of beams and shells and a linear modal basis.
However, similar techniques encounter severe difficulties when 
using general 3D elements, see e.g.\ \cite{vizzaccaro2020nonintrusive,givois2021dynamics,yichang2021comparison}.
For instance, the approach based on Implicit Condensation~\cite{hollkamp2008reducedorder}, as for instance proposed in \cite{frangi2019reduced,gobat2022ROM12ICEplus},
that has been tailored to resonant microdevices like gyroscopes or accelerometers in small transformations, fails when applied to micromirrors undergoing large rotations~\cite{opreni2021model}.
Similar difficulties are experienced with the quadratic manifold technique
built from modal derivatives~\cite{jain2017quadratic,vizzaccaro2020comparison,touze2021review}. 
Preliminary successful tests with the DPIM 
have been reported in~\cite{opreni2021model,vizza21high} using 
as actuation a fictitious modal forcing. 
The aim of this investigation is on the contrary to develop a reduction procedure 
for the Full Order Model accounting exactly for the piezo forcing. 
While we are not addressing here a fully coupled multiphysics problem,
this work is however intended as a first step in this direction and represents a major
achievement in itself as it provides the ideal simulation
tool for a whole class of high-impact technological applications.

In the developments reported herein, the focus is set on  piezo-MEMS actuators fabricated using Lead Zirconate Titanate (PZT), which is deposited in the form of a thin film sol-gel on bulk silicon. The resulting PZT is a random solid solution between PbTiO$_3$ and
PbZrO$_3$ and is widely used due to its excellent properties, 
such as high relative dielectric constants, high remnant polarisation and large piezoelectric coefficients.
In the applications targeted herein
typical PZT thin films have the thickness of few microns and are actuated with voltage biases in the order of tens of Volts generating oscillating electric fields with maximum values often larger than $10^7$\,V/m, which unavoidably exceed the linear range of piezoelectric materials. Therefore, models that predict the effects of nonlinearities are a necessary  tool for the design and simulation of this class of devices, and reduction methods accounting for these effects are a key tool at the design stage.
Although research is rapidly progressing, an ab-initio accurate
numerical computation of macroscopic hysteresis loops is still
beyond current capabilities as it is strongly influenced by microstructural defects arising from fabrication processes \cite{fedeli2019phase,fedeli2019defects}. For these
reasons we adopt herein a pragmatic approach based on direct
measurements of the polarisation through experiments, as proposed
in \cite{opreni2021modeling} where the time evolution of the average polarisation field within the piezoelectric material 
was measured using a standard Sawyer-Tower circuit for each value of the voltage applied.
The polarisation induces, through the converse piezoelectric effect, 
inelastic strains and stresses that actuate the device according to the fundamental elements of the Landau-Devonshire theory 
of ferroelectric materials~\cite{devonshire1954theory}.
These assumptions have been tested and validated with experiments on micromirrors
in \cite{frangi2020nonlinear} and following this approach a similar method has been implemented in 
commercial codes like Comsol Multiphysics\textsuperscript{\textregistered} \cite{comsol}.

The paper is organized as follows. 
After setting in Section~\ref{sec:equations} the equilibrium and constitutive equations
with the associated spectral properties, the DPIM theory is briefly recalled 
in Section~\ref{sec:dpim}, with specific emphasis on new features 
with respect to known developments.
Finally, in Section~\ref{sec:examples} two different examples
are discussed, starting with an academic benchmark 
on a doubly clamped beam.
Two MEMS micromirrors are then addressed, for which also experimental data 
are available. 
All the numerical examples are validated against an Harmonic Balance approach for the 
full-order model. Additional numerical results on the academic case of a cantilever beam are also collected in Appendix~\ref{app:cantilever} for the sake of completeness.

\section{Piezoelectric MEMS modelling}

This section is devoted to detailing the governing equations used to model a nonlinear structure actuated with piezoelectric effect considered in this study. The equations will be written starting from the strong form and then 
addressing the derivation of the weak form and the semi-discrete equations using 
a classical FE procedure. As a consequence of the piezoelectric forces, the static position of the structure at rest is modified. Hence a first necessary step consists in computing this static position, before analyzing the nonlinear vibrations.

\subsection{Governing equations}
\label{sec:equations}

Piezoelectric materials like PZT are characterized by a strong electrostrictive response.
The polarisation field induces inelastic strains $\eu^p$ 
so that the assumed Kirchhoff-like constitutive equation,
which holds under the classical assumption of large transformations and small strains,
reads:
\begin{equation}
\label{eq:electrostriction}
\Su = \cA:\eu-\Su^p,
\quad \text{with} \quad
 \Su^p := \cA:\eu^p,
\end{equation}
where $\Su$ is the second Piola-Kirchhoff stress tensor,
$\Su^p$ are inelastic stresses,
$\eu$ is the total Green Lagrange strain tensor
and $\cA$ is the fourth order elasticity tensor. 

According to the theory of electrostrictive materials introduced by Landau-Devonshire \cite{devonshire1954theory}, inelastic strains depend on the polarisation as follows:
\begin{equation}
\label{eq:epsp}
\eu^p = \cQ:(\pu\otimes\pu),
\qquad
e^p_{ij} = \sum_{k,\ell}\cQ_{ijk\ell}p_k p_{\ell}.
\end{equation}
Here $\pu$ denotes the polarisation back-rotated in the reference configuration and $\cQ$ is the electrostrictive coefficients tensor. 

It is worth stressing that in this work
the polarisation history at every point of the piezo patches 
is measured experimentally and is treated as a known periodic function of time.
In particular, we assume that the polarisation is not affected by the deformation of the device
which is a simplifying assumption holding anyway with very good accuracy for actuators undergoing moderate transformations.
The dynamic response of the structure is hence entirely defined by the conservation of linear momentum, the latter expressed in the reference configuration $B$:
\begin{equation}
 \rho \ddot{\uu} - \nabla\cdot\Pu=\rho\Bu, 
 \quad \text{in} \quad B\times t \in (t_0,T_e],\label{eq:strong_pde}
\end{equation}
with $\rho$ density, $\uu$ displacement field, $\Pu$ first Piola-Kirchhoff stress tensor,  $\nabla$ gradient operator, $\ddot{(\cdot)}$ second partial derivative with respect to time,  $\Bu$ body forces per unit mass. For the sake of simplicity we will neglect in the following every external action in the form of body or surface forces, as well as non-homogeneous kinematic boundary conditions, the only actuation being provided by the imposed polarisation through the electrostrictive effect. 
If needed, these could be easily included in the formulation. 
All quantities are defined in the reference configuration $B$ and over the time span from $t_0$ and $T_e$. Boundary and initial conditions for Eq.~\eqref{eq:strong_pde} are given as:

\begin{subequations}
\label{eq:bc_pde}
  \begin{align}
     \Pu\cdot\nuu=\zerou,     &  \qquad \text{on} \, \partial B^{N} \times t \in (t_0,T_e],
     \\
     \uu = \zerou,      &  \qquad \text{on} \, \partial B^{D} \times t \in (t_0,T_e],
     \\
     \uu = \uu_0,          &  \qquad \text{in} \, B\times t = t_0,
     \\
     \dot{\uu}=\dot{\uu}_0, & \qquad \text{in} \, B\times t = t_0,
  \end{align}
\end{subequations}

Equation \eqref{eq:strong_pde} is recast into a weak formulation upon introduction of a test function field $\wu$ which is defined over the space of admissible functions that vanish on the portion of the boundary where Dirichlet boundary conditions on $\uu$ are prescribed, i.e. $\cC(\zerou)$. Projecting the governing equation onto the test functions, the following weak formulation is obtained for the problem at hand:

\begin{align}
\label{eq:lin_mom_weak_full}
    & \int_{B} \rho\,\ddot{\uu}\cdot\wu\,\dd B + \int_{B} \eu:\cA:\delta\eu\,\dd B =  \int_{B_p} \Su^p:\delta\eu \,\dd B, \quad \forall\,\wu\in\cC(\zerou),
\end{align}

The inelastic stresses $\Su^p$ are integrated only over $B_p$,
defined as the collection of piezoelectric patches. 
Since $\Su^p$ are given periodic functions of time
providing the actuation mechanism, they are accordingly collected at the right-hand side of 
Eq.~\eqref{eq:lin_mom_weak_full}.
The Green-Lagrange strain tensor $\eu$ and its first variation $\delta\eu$ are 
defined as:

\begin{subequations}\label{eq:weak_pde}
  \begin{align}
    \eu =&\, \mathrm{sym}(\nabla\uu) + 
    \frac{1}{2}\nabla^{T}\uu\cdot\nabla\uu,  \\
    \delta\eu =&\, \mathrm{sym}(\nabla\wu) + \mathrm{sym}(\nabla^{T}\wu\cdot\nabla\uu).
  \end{align}
\end{subequations}

It is worth stressing that Eq.~\eqref{eq:lin_mom_weak_full} 
treats geometric nonlinearities exactly and its validity is only limited by the
small strain assumption of the Kirchhoff constitutive law~\eqref{eq:electrostriction}. 

All nonlinear terms are polynomial and the resulting expression upon explicit decomposition of such term is given as:

\begin{align}\label{eq:lin_mom_weak_full_explicit}
    & \int_{B} \rho\,\ddot{\uu}\cdot\wu\,\dd B + \int_{B} \sym(\nabla\uu):\cA:\sym(\nabla\wu)\,\dd B + \nonumber \\ 
    &  \int_{B} \sym(\nabla\uu):\cA:\sym(\nabla^{T}\wu\cdot \nabla\uu)\,\dd B + \frac{1}{2} \int_{B} \sym(\nabla\wu):\cA:\sym(\nabla^{T}\uu\cdot \nabla\uu)\,\dd B  + \nonumber \\
    & \frac{1}{2} \int_{B} \sym(\nabla^{T}\uu\cdot \nabla\uu):\cA:\sym(\nabla^{T}\wu\cdot \nabla\uu)\,\dd B = \nonumber \\
    & \int_{B} \Su^p:\sym(\nabla\wu) \,\dd B + \int_{B} \Su^p:\sym(\nabla^{T}\uu\cdot\nabla\wu) \,\dd B.
\end{align}

Equation \eqref{eq:lin_mom_weak_full_explicit} can be discretised using for instance the finite element method with nodal shape functions. Detailed derivation of the discretisation scheme of all quantities is reported in Appendix \ref{app:FEM_formulation}. Upon addition of linear damping the following system of time-dependent differential equations is derived:

\begin{equation}
\label{eq:lin_mom_full_discr0}
    \Mu \ddot{\Uu} + \Cu \dot{\Uu} + \Ku\Uu + \G(\Uu,\Uu) + \H(\Uu,\Uu,\Uu) = \Fu_P(t) + \Ku_P(t) \Uu,
\end{equation}

where $\Mu,\Cu,\Ku$ are respectively the mass, damping and stiffness matrices, $\G$ and $\H$ represent the quadratic and cubic nonlinearity tensors,  $\Fu_P$ represents the time-dependent piezoelectric force, while $\Ku_P$ stands for the time-dependent piezoelectric stiffness. 
The right-hand side of Eq.~\eqref{eq:lin_mom_full_discr0} stems from the discretisation of the contribution due to inelastic stresses induced by the polarisation.
Equation \eqref{eq:lin_mom_full_discr0} represents a non-autonomous dynamical system with nonlinear terms up to cubic order. Piezoelectric forcing yields stiffness terms that alter the eigenspectrum of the system and need to be properly treated during the reduction procedure.

\subsection{Computation of the static equilibrium position}
\label{sec:fixed}

The piezoelectric forces in Eq.~\eqref{eq:lin_mom_full_discr0} contain constant terms, which in turn create a new static position at rest for the structure. In order to compute the nonlinear vibration around this static position, one needs to split all time-dependent terms into their mean value over time and their time-dependent component. For the rest of the paper, it is assumed that the external excitation is periodic with period $T$. The decomposition of piezoelectric force and stiffness writes:

\begin{subequations}
  \begin{align}
    \Fu_P(t) =& \hFu_P^{(0)} + \hFu_P(t),\\
    \Ku_P(t) =& \Ku_P^{(0)} + \tKu_P(t),
  \end{align}
\end{subequations}

where the average value of the excitation is computed as:

\begin{subequations}
  \begin{align}
    \hFu_P^{(0)} =& \frac{1}{T} \int_{0}^{T} \Fu_P \, \dd t,\\
    \Ku_P^{(0)} =& \frac{1}{T} \int_{0}^{T} \Ku_P \, \dd t.
  \end{align}
\end{subequations}

As a result, one can clearly distinguish the autonomous terms from the non-autonomous ones, and collect the non-autonomous terms on the right-hand side of the equations. The final expression for the semi-discrete equations of motion is then expressed as:

\begin{align}\label{eq:lin_mom_full_discr}
     \Mu \ddot{\Uu} + \Cu \dot{\Uu} + (\Ku - \Ku_P^{(0)} )\Uu + \G(\Uu,\Uu) + \H(\Uu,\Uu,\Uu) - \hFu_P^{(0)} = \veps(\hFu_P + \tKu_P \Uu).
\end{align}

In this equation, a book-keeping parameter $\veps$ is added to make it compatible with the DPIM procedure for non-autonomous systems, as detailed in \cite{opreni22high}. 

The constant forcing terms in the left-hand side of Eq.~\eqref{eq:lin_mom_full_discr} needs to be first balanced in order to obtain the new static position of the structure at rest. Let us denote as $\Uu_0$ this deformed static configuration. Since one is interested in computing the nonlinear vibrations around this static position, the nodal displacement vector is expanded along:

\begin{equation}\label{eq:displ_expansion}
  \Uu(t) = \Uu_0 + \tUu(t),
\end{equation}

where $\tUu(t)$ is the time-dependent part of the displacement field. Plugging this ansatz into Eq.~\eqref{eq:lin_mom_full_discr}, one first obtains a static problem that needs to be solved for $\Uu_0$, which reads:

\begin{equation}\label{eq:static}
  (\Ku - \Ku_P^{(0)} )\Uu_0 + \G(\Uu_0,\Uu_0) + \H(\Uu_0,\Uu_0,\Uu_0) - \hFu_P^{(0)} = \zerou.
\end{equation}

The solution of this nonlinear static problem is obtained classically thanks to a Newton-Raphson procedure. In particular, an explicit analytical decomposition of the internal power into polynomial terms is not required for this step since the linearisation of the internal power can be performed directly through computation of the material stiffness and geometrical stiffness \cite{holzapfel2000nonlinear} as reported in Appendix \ref{app:static_analysis}.

Upon substitution of Eq.~\eqref{eq:displ_expansion} into Eq.~\eqref{eq:lin_mom_full_discr}, all static terms simplify and one finally obtains the dynamical equations governing the nonlinear vibrations around the static position as:

\begin{align}\label{eq:lin_mom_full_discrB}
    & \Mu \ddot{\tUu} + \Cu \dot{\tUu} + \left[ (\Ku - \Ku_P^{(0)} )\tUu + 2\G(\tUu,\Uu_0) + 3\H(\tUu,\Uu_0,\Uu_0)\right] + \nonumber \\
    & \left[ \G(\tUu,\tUu) + 3\H(\tUu,\tUu,\Uu_0)\right] + \H(\tUu,\tUu,\tUu) = \veps( \hFu_P + \tKu_P \Uu_0 + \tKu_P \tUu).
\end{align}

the following auxiliary quantities can now be introduced:

\begin{subequations}\label{eq:zitang00}
  \begin{align}
    \tKu \tUu =& (\Ku  - \Ku_P^{(0)} )\tUu + 2\G(\tUu,\Uu_0) + 3\H(\tUu,\Uu_0,\Uu_0),\label{eq:zitang00a}\\
    \tG(\tUu,\tUu) =&  \G(\tUu,\tUu) + 3\H(\tUu,\tUu,\Uu_0),\label{eq:zitang00b}
    \\
    \tFu_P = & \hFu_P + \tKu_P \Uu_0,\label{eq:zitang00c}
  \end{align}
\end{subequations}

that allows rewriting Eq.~\eqref{eq:lin_mom_full_discr} in more compact form as:

\begin{equation}\label{eq:lin_mom_full_compact}
    \Mu \ddot{\tUu} + \Cu \dot{\tUu} + \tKu \tUu + \tG(\tUu,\tUu) + \H(\tUu,\tUu,\tUu) = \veps(\tFu_P + \tKu_P \tUu).
\end{equation}

The reduction method based on the direct parametrisation of invariant manifold will be used to analyze the nonlinear vibrations displayed by Eq.~\eqref{eq:lin_mom_full_compact}. The autonomous part in the left-hand side contains quadratic and cubic nonlinearities, fitting with the generic framework given in~\cite{vizza21high}. The non-autonomous terms in the right-hand side differ from those treated in~\cite{opreni22high} by a term that is linear with respect to the displacement and needs a dedicated development. In particular, this term has a direct influence on the eigenfrequencies of the system, nevertheless it is treated as a small excitation term in the DPIM procedure, which is computed for the eigenvalues corresponding to the autonomous terms on the left-hand side. This assumption will be shown to give effective results, provided that the perturbation on the eigenfrequencies is small. This will be verified quantitatively in the numerical examples.

All the spectral properties of the mechanical system detailed in \cite{vizza21high,opreni22high} 
also apply here, the only difference being that the spectrum of the system is computed at the new fixed point, as also recently processed for rotating structures in~\cite{Martin:rotation}. The resulting orthogonality properties of the system are expressed with respect to the tangent stiffness $\tKu$ since it corresponds to the linear stiffness operator at the configuration where the system is linearised. In order to properly derive the ROM,  Eq.~\eqref{eq:lin_mom_full_compact} is first rewritten as a first-order dynamical system

\begin{subequations}\label{eq:dpim_start}
  \begin{align}
    & \Mu \dot{\tVu} + \Cu \tVu + \tKu \tUu + \tG(\tUu,\tUu) + \H(\tUu,\tUu,\tUu) = \veps(\tFu_P + \tKu_P \tUu),\\
    &\Mu\dot{\tUu} = \Mu \tVu,
  \end{align}
\end{subequations}
where the velocity $\tVu$ associated to $\tUu$ is identical to the physical velocity since $\dot{\Uu}_0=\zerou$.

\subsection{Spectral properties}

Let us introduce the eigenfunctions $\bPhi$ end eigenfrequencies $\omega_j$ as the solution of the following problem:

\begin{equation}
    \left( -\omega_j^2\Mu+\tKu \right) \bPhi_j = \zerou.
\end{equation}

The tangent stiffness operator is symmetric, hence eigenfunctions and eigenvalues are real valued. Hereafter, we assume that the tangent-operator remains positive-definite (instabilities like  buckling are not considered). The left and right  eigenfunctions, denoted as  $\bX_j$ and $\bY_j$, together with  the eigenvalues $\Lambda_j$ of the first-order problem, are defined as the solutions of the following systems:

\begin{equation}
\left(
\Lambda_s
\begin{bmatrix}
\Mu & \zerou
\\
\zerou & \Mu
\end{bmatrix}
+
\begin{bmatrix}
\Cu & \tKu
\\
-\Mu & \zerou
\end{bmatrix}
\right)
\bY_s
=
\zerou, \qquad
\bX_s^\text{T}
\left(
\Lambda_s
\begin{bmatrix}
\Mu & \zerou
\\
\zerou & \Mu
\end{bmatrix}
+
\begin{bmatrix}
\Cu & \tKu
\\
-\Mu & \zerou
\end{bmatrix}
\right)
=\zerou,
\end{equation}

where properties and explicit definitions for right and left eigenvalues are detailed in \cite{vizza21high,opreni22high}. Let us now define the set of master modes on which the reduction method will rely. Assuming that $n$ master modes are selected for the ROM, 
we introduce the following matrices related to the master linear subspace:

\begin{subequations}
    \begin{align}
\Xu =&
    \begin{bmatrix}
    \bX_{m_1} &\bX_{m_2} &\ldots &\bX_{m_n}  
    &\bar{\bX}_{m_1} &\bar{\bX}_{m_2} &\ldots& \bar{\bX}_{m_n} 
    \end{bmatrix},
\\
\Yu = &
    \begin{bmatrix}
    \bY_{m_1} &\bY_{m_2} &\ldots &\bY_{m_n}  
    &\bar{\bY}_{m_1} &\bar{\bY}_{m_2} &\ldots& \bar{\bY}_{m_n} 
    \end{bmatrix},
\\
\lambdau = &
\text{diag}[\Lambda_{m_1},\;\Lambda_{m_2},\;\ldots,\;\Lambda_{m_n}\;\bar{\Lambda}_{m_1},\;\bar{\Lambda}_{m_2},\;\ldots,\;\bar{\Lambda}_{m_n}],
\label{eq:master_eigs_HO}
\\
\phiu = &
\begin{bmatrix}
\bPhi_{m_1} &\bPhi_{m_2} &\ldots &\bPhi_{m_n}  &\bPhi_{m_1} &\bPhi_{m_2} &\ldots& \bPhi_{m_n} 
\end{bmatrix},
\label{eq:master_eigv_HO}
\end{align}
\end{subequations}

with $\Xu$ and $\Yu$ the $2N\times 2n$ matrices of left and right master eigenvectors, $\lambdau$ the $2n\times2 n$ matrix of complex master eigenvalues, and $\phiu$ the $N\times 2n$ matrix of master modes. It follows that the newly introduced matrices can be also written as:

\begin{subequations}\label{eq:master_leftright_HO}
\begin{align}
    \Xu =& \begin{bmatrix}
    \Xu_1 &\Xu_2 &\ldots &\Xu_n  
    &\bar{\Xu}_1 &\bar{\Xu}_2 &\ldots& \bar{\Xu}_n 
    \end{bmatrix},\\
    \Yu =& \begin{bmatrix}
    \Yu_1 &\Yu_2 &\ldots &\Yu_n  
    &\bar{\Yu}_1 &\bar{\Yu}_2 &\ldots& \bar{\Yu}_n 
    \end{bmatrix},\\
    \lambdau =&
    \text{diag}[\lambda_1,\,\lambda_2,\,\ldots,\lambda_n,\,\bar{\lambda}_1,\,\bar{\lambda}_2,\,\ldots,\,\bar{\lambda}_n],\\
    \phiu =& \begin{bmatrix}
    \phiu_1 &\phiu_2 &\ldots &\phiu_n  &\phiu_1 &\phiu_2 &\ldots& \phiu_n 
    \end{bmatrix},
\end{align}
\end{subequations}

which makes clear that in the sorting of the master quantities, the $j$-th master index corresponds to the $m_j$-th index in the sorting of the original system and the $(j+n)$-th to the $(m_j+N)$-th, with $j\in(1,n)$.

\section{Reduced-order modelling strategy}

In this Section, the main equations needed to derive the ROM, are detailed. The method relies on the computation of a nonlinear mapping and the reduced dynamics, both being computed with polynomial expansions at arbitrary order. The complete derivation of the method is developed in~\cite{vizza21high,opreni22high} and only the main points are here briefly recalled by underlining the additional steps needed to tackle the peculiarity of the problem at hand.

\subsection{Direct parametrisation of invariant manifold}
\label{sec:dpim}

The direct parametrisation method relies on the introduction of nonlinear mappings relating initial nodal displacement and velocity vectors $\tUu$ and $\tVu$ in physical space, to newly introduced {\em normal coordinate} $\zu$, describing the motions on the invariant manifold associated to the selected master coordinates. Since the problem at hand is non-autonomous, one can use the nonlinear mappings introduced in~\cite{opreni22high}, which read:

\begin{subequations}
\label{eq:zeNLmaps}
  \begin{align}
    \tUu = \WUfunsim(\zu,\Omegau,t) = \WUfun + \veps \hWUfunsim(\zu,\Omegau,t) + \cO(\veps^2),\\
    \tVu = \WVfunsim(\zu,\Omegau,t) = \WVfun + \veps \hWVfunsim(\zu,\Omegau,t) + \cO(\veps^2).
  \end{align}
\end{subequations}

Importantly, the change of coordinates are composed of two different terms, the first one at order $\varepsilon^0$ being related to the autonomous problem, while the second one at order $\varepsilon^1$ is concerned with the non-autonomous forcing terms. For the last term, one can note the dependence upon time $t$ and upon the excitation frequencies $\Omegau$. The reduced dynamics along the embedding can be also be decomposed as:

\begin{equation}\label{eq:reduceddyn00}
    \dot{\zu} = \ffunsim(\zu,\Omegau,t) = \ffun + \veps\hffunsim(\zu,\Omegau,t)+\cO(\veps^2),
\end{equation}

where the corresponding splitting pertaining to autonomous and non-autonomous terms is used accordingly. As a consequence, the reduced dynamics depends explicitly on time and on parameters as the excitation frequencies $\Omegau$ as well. The reduced dynamics and the mappings depend also on implicit parameters such as geometry and material parameters, which affect the resulting values of the coefficients of the vector field $\ffunsim(\zu,\Omegau,t)$. 

In order to compute the unknown mappings and reduced dynamics, the starting point consists in deriving the {\em invariance equation}~\cite{cabre2005overview,haro2016parametrisation}, which is found by eliminating time. This equation embeds the invariance property of the searched manifolds where the reduced dynamics will lie. It is simply found by differentiating Eq.~\eqref{eq:zeNLmaps} with respect to time, substitute into~\eqref{eq:dpim_start} and use~\eqref{eq:reduceddyn00}. In the non-autonomous case, terms of like-powers of $\varepsilon$ are also collected, such that two different invariance equations are written~\cite{opreni22high}. At order $\varepsilon^0$, the invariance equation for the autonomous terms writes:

\begin{subequations}
\label{eq:invariance_zero}
  \begin{align}
  	& \Mu\nabla_{\zu}\WVfun\ffun + \Cu\WVfun + \tKu\WUfun + {\tG}(\WUfun,\WUfun) + \H(\WUfun,\WUfun,\WUfun) = \zerou, \\
  	& \Mu\nabla_{\zu}\WUfun\ffun - \Mu\WVfun = \zerou.
  \end{align}
\end{subequations}

One can note in particular that Eq.~\eqref{eq:invariance_zero} is identical to the one handled in~\cite{vizza21high}, hence the same algorithm can be adopted.  A different result is observed for the  $\veps^1$-invariance equation due to the extra term in the right-hand side of~\eqref{eq:dpim_start} which is proportional to the displacement, yielding:

\begin{subequations}
\label{eq:invariance_one}
  \begin{align}
  	& \Mu\left(\dot{\hWVfunsim} (\zu,\Omegau,t) + \nabla_{\zu}\WVfun\hffunsim(\zu,\Omegau,t) + \nabla_{\zu}\hWVfunsim(\zu,\Omegau,t)\ffun  \right) + \Cu\hWVfunsim(\zu,\Omegau,t) + \tKu\hWUfunsim(\zu,\Omegau,t) \nonumber \\
  	& + 2\tG(\WUfun,\hWUfunsim(\zu,\Omegau,t)) + 3 \H(\WUfun,\WUfun,\hWUfunsim(\zu,\Omegau,t)) = \tilde{\Fu}_P + {\tKu}_P\WUfun,\label{eq:invariance_one_a}\\
  	& \Mu\left(\dot{\hWUfunsim} (\zu,\Omegau,t) + \nabla_{\zu}\WUfun\hffunsim(\zu,\Omegau,t) + \nabla_{\zu}\hWUfunsim(\zu,\Omegau,t)\ffun  \right) - \Mu \hWVfunsim(\zu,\Omegau,t) = \zerou.
  \end{align}
\end{subequations}

Contrary to the developments reported in~\cite{opreni22high}, an additional term ${\tKu}_P\WUfun$ is present in the right-hand side of Eq.~\eqref{eq:invariance_one_a}, as a consequence of the time-dependent piezoelectric stiffness. Equation~\eqref{eq:invariance_one} represents a system of linear differential equations and can be solved using Fourier analysis. To this aim, let us decompose the external excitation terms in its Fourier components:

\begin{subequations}
\label{eq:decompforce}
  \begin{align}
  	\tilde{\Fu}_P = \sum_{j=1}^{2\hn} \tilde{\Fu}_{P_j} e^{\hlm_{j}t},\\
  	\tilde{\Ku}_P = \sum_{j=1}^{2\hn} \tilde{\Ku}_{P_j} e^{\hlm_{j}t},
  \end{align}
\end{subequations}

where $\hlm_{j}$ represents either $+\iu\Omega_j$ or $-\iu\Omega_j$. In practice, if no internal resonances are experienced by the structure, only the first harmonic component is necessary to derive an accurate reduced model. Equations~\eqref{eq:invariance_one} is linear with respect to the non-autonomous mappings and reduced dynamics, hence the resulting quantities will have the same frequency content:

\begin{align}
    \hWUfunsim(\zu,\Omegau,t) = \sum_{j=1}^{2\hn} \hWUfunsim_j(\zu) e^{\hlm_{j}t}, \quad \hWVfunsim(\zu,\Omegau,t) = \sum_{j=1}^{2\hn} \hWVfunsim_j(\zu) e^{\hlm_{j}t}, \quad \hffunsim(\zu,\Omegau,t) = \sum_{j=1}^{2\hn} \hffunsim_j(\zu) e^{\hlm_{j}t},
\end{align}

where $\hWVfunsim_{j}(\zu)$, $\hWUfunsim_{j}(\zu)$ and  $\hffunsim_{j}(\zu)$ are coefficients of the Fourier expansion. We highlight that the coefficients are neither function of time nor of the frequency. Upon substitution of these expansions in Eq.~\eqref{eq:invariance_one} we can then project the system onto Fourier basis, hence providing the following representation for the $\veps^{1}$-invariance equation:

\begin{subequations}
\label{eq:invariance_one_proj}
  \begin{align}
  	& \forall \, j= 1,...,2\hn,\nonumber \\
  	& \hlm_{j} \Mu \hWVfunsim_{j}(\zu) + \Mu \nabla_{\zu} \hWVfunsim_{j}(\zu) \ffun + \Mu \nabla_{\zu} \WVfun \hffunsim_{j}(\zu) + \Cu\hWVfunsim_{j}(\zu) \nonumber \\
  	& + \tKu \hWU_{j}(\zu) + 2{\tG}(\hWUfunsim_{j}(\zu),\WUfun) + 3\H(\hWUfunsim_{j}(\zu),\WUfun,\WUfun) =  \tilde{\Fu}_{P_j} + \tilde{\Ku}_{P_j}\WUfun, \\
  	& \hlm_{j} \Mu \hWUfunsim_{j}(\zu) + \Mu \nabla_{\zu} \hWUfunsim_{j}(\zu) \ffun + \Mu \nabla_{\zu} \WUfun \hffunsim_{j}(\zu) - \Mu \hWVfunsim_{j}(\zu) = \zerou,
  \end{align}
\end{subequations}

which can be solved recursively to compute the Fourier coefficients associated to mappings $\hWVfunsim_{j}(\zu)$, $\hWUfunsim_{j}(\zu)$ and to the reduced dynamics $\hffunsim_{j}(\zu)$. One can remark that, following the projection onto Fourier basis of the different harmonic components, the invariance property is correctly recovered for the $\varepsilon^1$ problem describing the non-autonomous part, which does not depend on time anymore. An important aspect of novelty as compared to past developments is that the $\veps^{1}$-invariance equation in presence of a time-dependent modulation of the stiffness features an additional term on its right-hand side: $\tilde{\Ku}_{P_j}\WUfun$. This last term is provided by the time-dependent part of the piezoelectric stiffness. 
We finally remark that the computed mappings and reduced dynamics coefficients show a dependence on the value of the excitation frequency since Eq.~\eqref{eq:invariance_one_proj} is obtained by projecting the system onto a proper Fourier basis. This dependence is mild and it can be neglected if the system is excited at resonance, as also highlighted in \cite{opreni22high}. 
Indeed, in the present developments we will show that highly accurate reduced models can be obtained by parametrising the system for a single excitation frequency value and then exploit the reduced model to compute the entire Frequency Response Curve (FRC) of the system.

\subsection{Solution scheme}

The solutions to Eqs. \eqref{eq:invariance_zero} and \eqref{eq:invariance_one_proj} are found expressing the unknown $\zu$ with arbitrary order polynomial expansions. This choice is mainly guided by the fact that recursive solutions, order by order, are possible and offer an accurate solution scheme. The invariance equations are indeed rewritten order by order, leading to the so-called {\em homological equations}. For the autonomous mappings, the expansions are searched for according to:

\begin{subequations}
  \begin{align}
      \WUfun =&\, \phiu \zu + \sum_{p=2}^{o}[ \WUfun ]_{p},\\
      \WVfun =&\, \phiu \lambdau \zu + \sum_{p=2}^{o}[ \WVfun ]_{p},
  \end{align}
\end{subequations}

where $o$ stands for the maximum order of the selected expansion, and the shortcut notation $[\cdot]_{p}$ indicates a generic term of order $p$. Importantly, the first linear term of these expansions underline that the mappings are identity-tangent to the master eigenmodes, for small vibration amplitudes, which is in line with the notion of a Nonlinear normal mode (NNM) defined as an invariant manifold tangent at origin to the master vibration modes. For the Fourier coefficients of the non-autonomous mappings, the following expansions are introduced, where the dependence on time is not present and the dependence on the excitation frequency is implicit, hence neither of them are reported:

\begin{subequations}
  \begin{align}
      \hWUfunsim_{j}(\zu) =&\, \hWUfun = \sum_{p=0}^{q}[ \hWUfun ]_{p},\\
      \hWVfunsim_{j}(\zu) =&\, \hWVfun = \sum_{p=0}^{q}[ \hWVfun ]_{p}.
  \end{align}
\end{subequations}

In these equations, the subscript $j$, referring to one of the excitation frequency $\Omega_j$ has been dropped for simplicity, since the same independent computation needs to be tackled for each of the driving frequencies. The order of the development of the non-autonomous mappings is $q<o$~\cite{opreni22high}. Note also that the series expansion starts at order~$0$ for the non-autonomous part,  corresponding to the simplest approximation that can be used to deal with the forcing, see discussions in {\em e.g.}~\cite{touze2006nonlinear,breunung2018explicit,ponsioen2019analytic,opreni22high}. Together with the mappings, arbitrary order Taylor expansions are also used to represent the reduced dynamics along the embedding: 

\begin{subequations}
  \begin{align}
      \ffun  = &\, \lambdau \zu + \sum_{p=2}^{o}[\ffun]_{p},\\
      \hffunsim_{j}(\zu) = \,& \hffun =  \sum_{p=0}^{q}[\hffun]_{p}.
  \end{align}
\end{subequations}

The next step consists in plugging these ansatz into the two invariance equations~\eqref{eq:invariance_zero} and \eqref{eq:invariance_one_proj}, for $\varepsilon^0$ and $\varepsilon^1$ orders. The remainder of this detailed calculation is reported in Appendix~\ref{app:ROMdetail} for the sake of brevity, since it follows the main guidelines provided in~\cite{vizza21high,opreni22high}. A special emphasis is put on processing the new terms appearing in the present problem. Identification of like-power terms leads to write an order-$p$ {\em homological equation}, which is also specifically written at the level of an arbitrary monomial present in the polynomial expansions, and in such a way that  direct computation from the physical space is possible. A key feature of these homological equations lies in the fact that they are underdetermined and ill-conditioned. The ill-conditioning is intimately linked to the notion of resonance~\cite{haro2016parametrisation,touze2014normal,haller2016nnm}. In the context of the direct calculation, this issue is solved by using a bordering technique~\cite{vizza21high}. The under-determinacy leads to an infinity of solutions, framed by the two most opposite styles given respectively by the graph style and the normal form style, see {\em e.g.}~\cite{haro2016parametrisation,vizza21high,opreni22high,touze2021review,jain2021how,StoychevRomer} for discussions. In the context of the present work, simulation results are presented using the complex normal form style.

The outcome of the whole reduction process can be summarized as follows. First, the addition of the piezoelectric force to the governing equations creates a static equilibrium position that has to be computed. Then the nonlinear vibrations around this deformed state are computed thanks to the DPIM, using an arbitrary order expansion of order $o$ for the autonomous part, and order $q < o$ for the non-autonomous part. In the remainder of the paper, the notation DPIM-${\mathcal O}(o,q)$ will be used to denote the orders retained. The reduced dynamics is given by a polynomial expansion of order $o$, and the complex normal form is used. The ROM equations are solved with numerical continuation, and the equations are realified following the procedure explained in~\cite{vizza21high,opreni22high}. In all the examples shown below, reduction to a single NNM is also used by keeping a single master coordinate.

\section{Applications}
\label{sec:examples}

In the present Section, the DPIM is first applied to
an academic example represented by a doubly clamped beam and 
then to real MEMS micromirrors 
developed by STMicroelectronics\texttrademark.
Additional results concerning a cantilever beam are collected in
Appendix \ref{app:cantilever}.
In what follows the results obtained with the proposed DPIM are  validated against
a full-order large-scale Harmonic Balance approach (HBFEM) that has been
developed and discussed in \cite{opreni2021analysis}. 
The HBFEM is considered the reference tool for these applications as far 
as the accuracy is concerned, even if its high computational cost hinders its applicability 
for the design of new devices and their optimization. 
On the contrary, the response of the reduced order models
provided by the DPIM is computed 
applying numerical continuation of periodic orbits with the MATCONT package \cite{dhooge2003matcont}.

The bulk structure of the devices analysed is
made either of polysilicon or of single crystal silicon
and their properties will be defined case by case.
On the contrary,
specific care has to be devoted to the treatment of the polarisation field in the piezo patches.

\subsection{Polarisation field}

The application of the proposed formulation 
requires the knowledge of the
periodic polarisation field at every point of the piezo patches and for every
time instant of the period. In what follows we briefly explain the pragmatic though accurate approach that is followed in this investigation in order to compute
the right-hand side in Eq.~\eqref{eq:lin_mom_full_discr0}. 

\begin{figure}[ht]
    \centering
    \includegraphics[width = .99\linewidth]{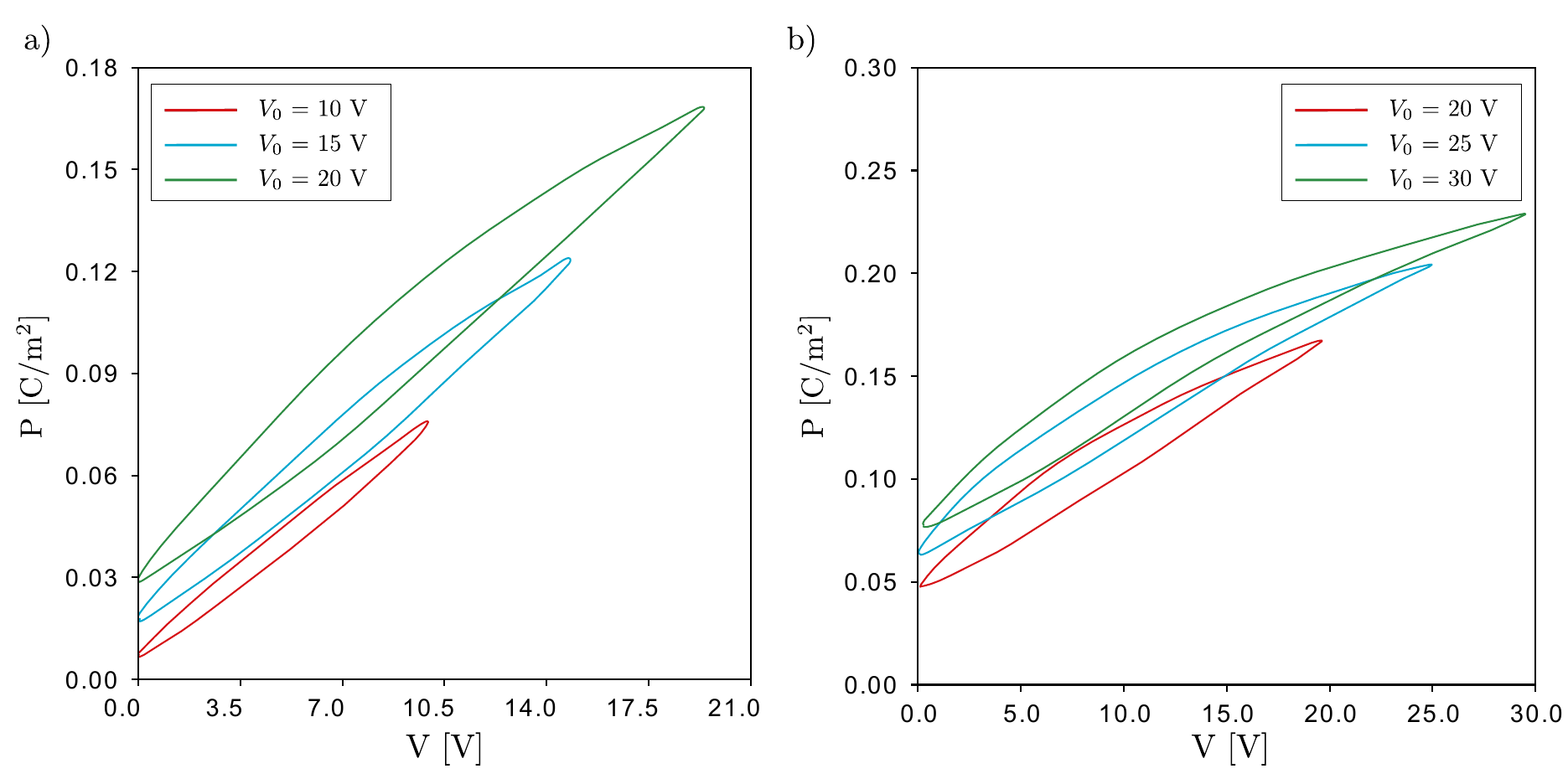}
    \caption{Polarisation curves measured for different potential values on two MEMS micromirrors that are fully described in Section~\ref{sec:mirrors} and are labelled as (a): Mirror A (a) and (b): Mirror B;
    see also Fig.~\ref{fig:mirrors_img} for their representations. The loops run in the counter-clockwise direction. These curves are generated by unipolar (i.e.\ always positive) voltage histories to avoid continuous switching of the polarisation domains in the piezo patch.}
    \label{fig:polarisation_curves}
\end{figure}

\begin{figure}[ht]
    \centering
    \includegraphics[width = .99\linewidth]{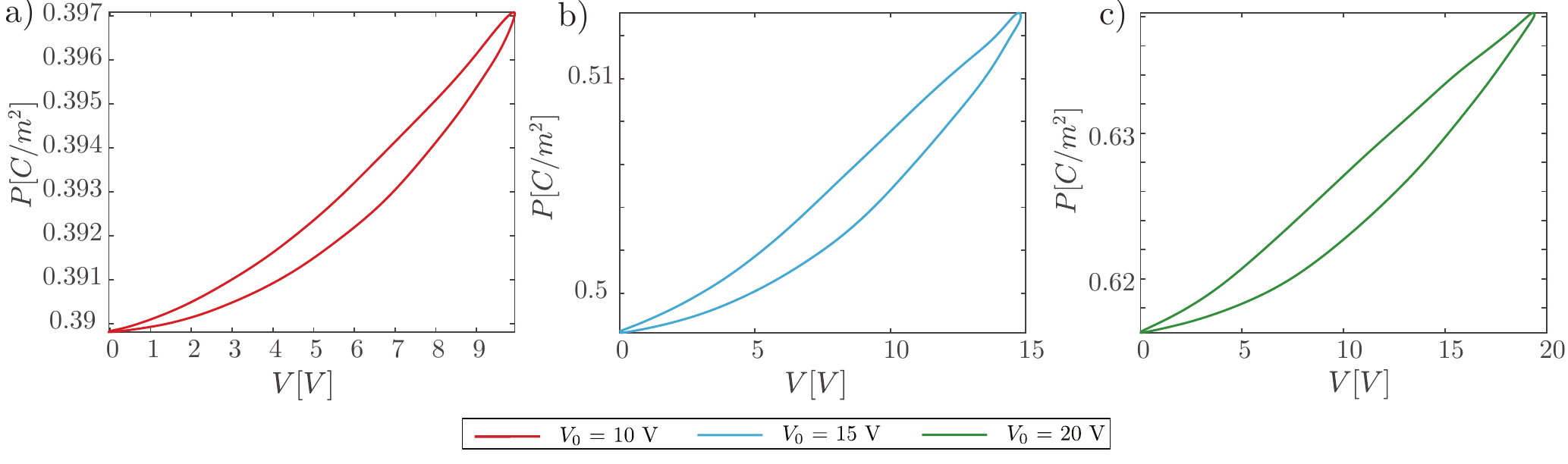}
    \caption{Polarisation curves modified to enhance the frequency shift effect for different potential values. The three figures refer to 10,15, and 20V respectively.}
    \label{fig:polarisation_curves_mod}
\end{figure}

During production, the piezo patches are deposited on top of the bulk of the 
solid on a plane of unit normal $\eu_3$. 
Their thickness $h$ is very small as compared to the in-plane
dimensions and are enclosed within upper and lower electrodes on which given voltage histories are imposed, generating a voltage bias $V(t)$.  As a consequence, the electric field in the piezo patches is almost 
exactly aligned with the $x_3$ direction and has intensity $V/h$.
Even if this does not imply that the polarisation itself is uniform within the films, 
the total force can be computed with very good accuracy 
by assuming that the field can be homogenised within 
a patch and that it has the form:
\begin{equation}
	\pu = P(t)\eu_3,
\end{equation}
with $P$ a scalar value. 
All these assumptions have been validated through an extensive experimental campaign
in \cite{opreni2021modeling,opreni2020piezoelectric} 
where it has also been evidenced that, for the working frequencies of interest, 
the polarisation history is frequency-independent, thus further reducing 
the experimental overhead required to initiate the simulation.
As an example, Fig.~\ref{fig:polarisation_curves} reports the
measured polarisation history $P(t)$  for two different types of PZT mixtures
utilized in the micromirrors discussed in Section \ref{sec:mirrors}.
The curves, which run in the counter-clockwise direction, highlight the typical pattern of polarisation in piezoelectric materials subjected to strong electric field values, i.e.\ strong hysteretic behaviour with lack of reversibility
and voltage dependence.
An important remark is that the polarisation measurements were performed in 
unipolar conditions, i.e.\  imposing a positive voltage bias of the type

\begin{align}
\label{eq:vbias}
V(t) = \frac{V_0}{2}(1+\cos{(\Omega t)}),
\end{align}

to avoid fatigue phenomena in the piezoelectric film
associated to continuous polarisation switching.
However, the proposed procedure could be applied to arbitrary polarisation cycles.
In all the applications discussed in the following sections, the PZT patches are divided in two sets labeled PZT-A and PZT-B, respectively. In the figures showing the geometries and the locations of these patches, PZT-A  are displayed with red colour patches while PZT-B with yellow patches, see {\em e.g.} Figs.~\ref{fig:ccbeam_geo} and \ref{fig:mirrors_img}.
These are subjected to unipolar potential histories of the type of
Eq.\eqref{eq:vbias} with phase shift:

\begin{equation}
\label{eq:actuation_law}
  V_{A} = \frac{V_0}{2}(1+\cos{(\Omega t)}), \quad V_{B} = \frac{V_0}{2}(1-\cos{(\Omega t)}),
\end{equation}

to maximise actuation.
In order to highlight some important effects of the piezoelectric actuating forces, we will also consider modified fictitious versions of the polarisation curves of Fig.~\ref{fig:polarisation_curves}. 
These new curves are tailored to enhance the frequency shift introduced by the piezoelectric forcing, which cannot be modelled with simplified techniques. 
The resulting hysteresis loops, corresponding to 10, 15 and 20V, are plotted in Fig.\ref{fig:polarisation_curves_mod} (a)-(c).

As a result, 
the piezoelectric strains and stresses defined in Eqs.\eqref{eq:electrostriction}-\eqref{eq:epsp} can be expressed everywhere in the PZT in an explicit manner 
and in terms of the known polarisation history.
Indeed, assuming transverse isotropy for the electrostrictive response of the thin film, the only non-zero components of the inelastic strains are:

\begin{gather}
\label{eq:epsP}
e^p_{11} =  e^p_{22} = \cQ_{1133}P^{2},
\quad 
e^p_{33} = \cQ_{3333}P^{2},
\end{gather}

where electrostrictive coefficients are taken from~\cite{haun1987thermodynamic} for a mole fraction
$x=0.5$ of PbTiO$_3$.

Finally, the inelastic stress to be inserted in
Eq.~\eqref{eq:lin_mom_weak_full_explicit} follows from the Kirchhof-Love assumptions as:

\begin{align}
    &S^p_{11} = S^p_{22} = \cA_{1111}e^p_{11}+\cA_{1122}e^p_{22}+\cA_{1133}e^p_{33} = \alpha_1 P^{2},
\nonumber \\
    &S^p_{33} = \cA_{3333}e^p_{33}+2\cA_{1133}e^p_{11} = \alpha_3 P^{2}.
\label{eq:sigmaP}
\end{align}
 
Very limited data are available concerning the elastic constants in
Eq.~\eqref{eq:sigmaP} for the type of PZT employed herein, but a good agreement with experiments
can be achieved 
with an assumption of isotropic material behaviour. 
Starting from Eq.~\eqref{eq:sigmaP} and applying standard FEM discretisation techniques, 
vector $F_p$ and matrix $K_p$ in
Eq.~\eqref{eq:lin_mom_full_discr0} can be readily computed for every time instant $t$,
and next decomposed in harmonic contributions according to Eq.~\eqref{eq:decompforce} 
in order to feed the reduction procedure detailed in Section \ref{sec:dpim}.

\subsection{Doubly clamped beam}
\label{sec:CCbeam}

\begin{figure}[htb!]
    \centering
    \includegraphics[width = .8\linewidth]{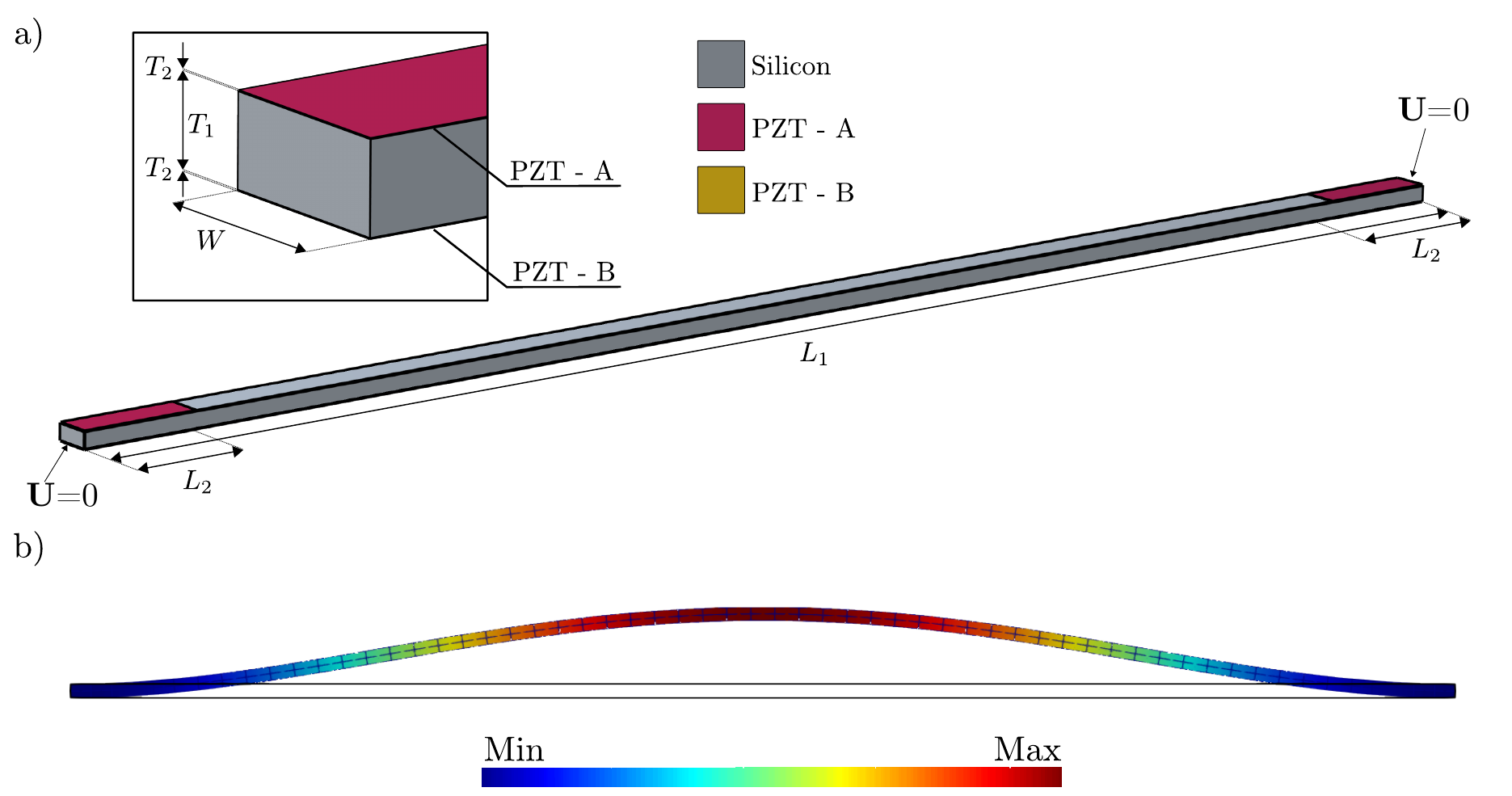}
    \caption{(a) Geometry of the tested clamped-clamped beam. 
    $T_1 = 1\,\mu$m, $T_2 = 0.01\,\mu$m, $L_1 = 100\,\mu$m, $L_2 = 7.5\,\mu$m, 
    $W = 2.2\,\mu$m. 
    (b) Shape of the first bending mode $\bPhi_B$. Dirichlet homogeneous boundary conditions  are imposed at the beam extremities.}
    \label{fig:ccbeam_geo}
\end{figure}

The first validation is performed on the clamped-clamped beam illustrated in
Figure~\ref{fig:ccbeam_geo}a, where $L_1 = 100\,\mu$m is the total length of the beam, 
$L_2 = 7.5\,\mu$m is the length of the four piezo patches.
$T_1 = 1\,\mu$m denotes the thickness of the silicon body of the beam, while 
$T_2 = 0.01\,\mu$m is the piezo thickness.

The material properties are detailed in Table \ref{tab:parameters_cc}.
In this academic test, both PZT and silicon have an isotropic mechanical behaviour.
Patches of type A are placed on the upper surface, while 
patches of type B, not visible in the Figure, are deposited on the lower surface
of the beam.
The two patches are actuated according to Eq.\eqref{eq:actuation_law}
and set in resonant motion the first bending mode illustrated 
in Figure~\ref{fig:ccbeam_geo}b. 
A mass-proportional Rayleigh damping model is considered with a quality factor $Q=100$.
The eigenfrequency of the first bending mode is  $\omega_0=5.399$\,rad/$\mu$s.

\begin{table}[h]
	\begin{center}
		\begin{tabular}{|ccc|}
              \hline
              \multicolumn{3}{|c|}{PZT}\\
			\hline
			Coefficient & Value & Unit\\
			\hline
			$\cQ_{3333}$ 		& 0.097 	& m$^4$/C$^2$ 		\\
			$\cQ_{1133}$ 		& -0.046 	& m$^4$/C$^2$		\\
			$E$ 		& 70000	& MPa 		\\
			$\nu$ 		& 0.33 	& -		\\
			\hline
                \hline
                \multicolumn{3}{|c|}{Silicon}\\
                \hline                
			Coefficient & Value & Unit\\
			\hline
			$E$ 		& 160000	& MPa 		\\
			$\nu$ 		& 0.22 	& -		\\
			\hline
		\end{tabular}
	\end{center}
	\caption{Constitutive parameters of PZT and silicon for the clamped-clamped and  cantilever beam examples. Isotropic mechanical behaviour for both materials is assumed. }
	\label{tab:parameters_cc}
\end{table}

The aim of this application is twofold. First, demonstrate the accuracy of the proposed DPIM formulation and, second, discuss important effects of the piezoelectric forcing terms that cannot be accounted for by earlier simplified formulations \cite{opreni22high}. To demonstrate the former we will consider the polarisation curves reported in 
Fig.~\ref{fig:polarisation_curves}a). The latter will be highlighted using the modified polarisation curves illustrated in Fig.~\ref{fig:polarisation_curves_mod}.

\begin{figure}[htb!]
    \centering
    \includegraphics[width = .7\linewidth]{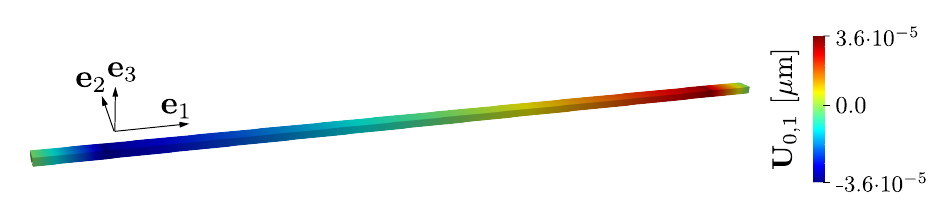}
    \caption{$\eu_1$ displacement field component associated to the fixed point of the clamped-clamped beam for $V_0$ = 20 V. The static load introduces an axial pre-stress which shifts the eigenfrequency.}
    \label{fig:clamped_fixed}
\end{figure}

The imposed polarisation oscillates around an average value that generates a static force component. This effect can be appreciated by inspecting Fig.~\ref{fig:clamped_fixed} which reports the static displacement field along the beam axis in the new fixed point. 
This induces an axial stress which stiffens the beam and shifts the eigenfrequency upwards.

\begin{figure}[htb!]
    \centering
    \includegraphics[width = .99\linewidth]{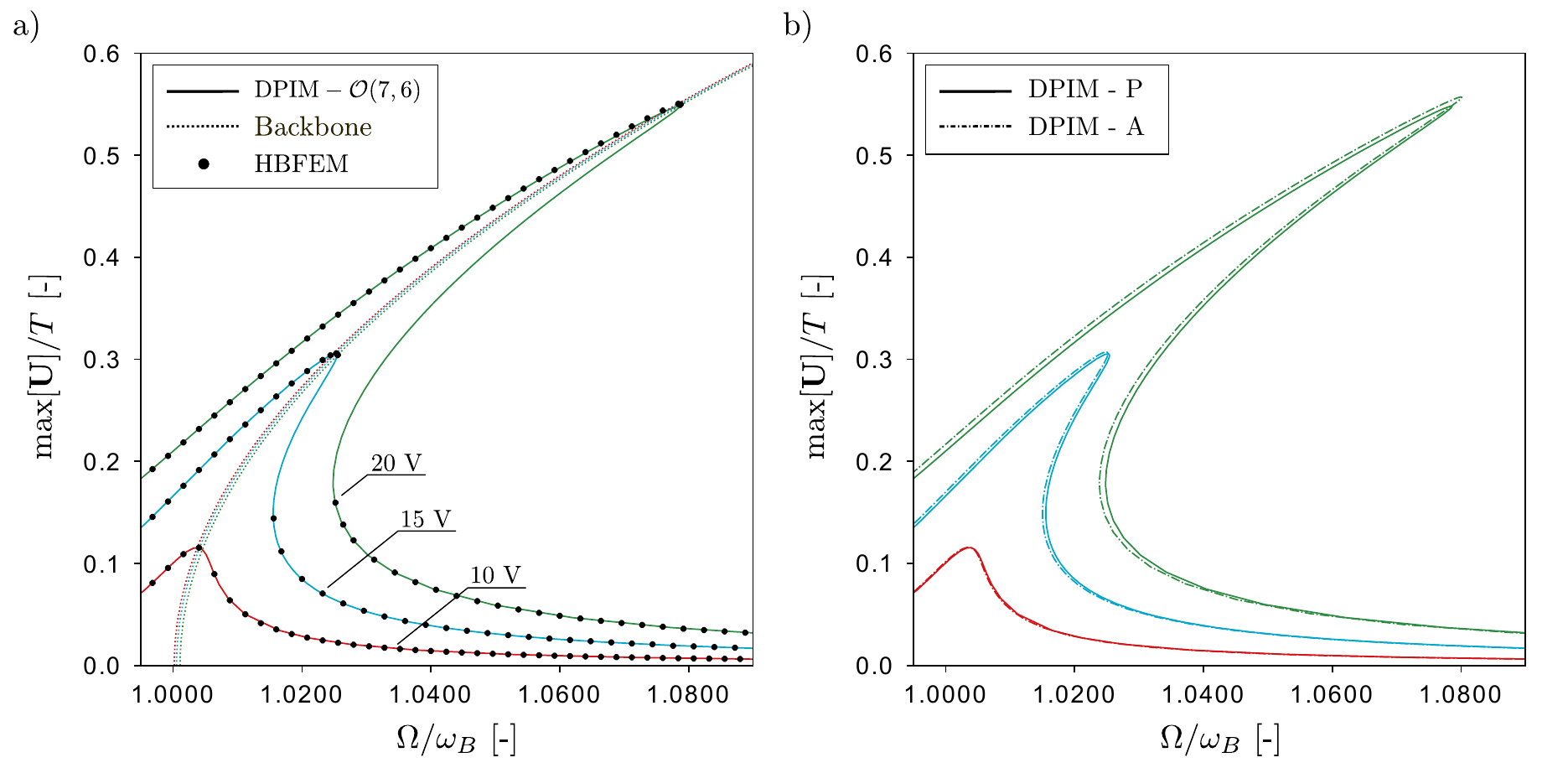}
    \caption{(a) Comparison between full order HBFEM simulations and DPIM reduced-order model computed for 10, 15, and 20~V. Backbone curves are reported to highlight the effect of the mean value of the piezoelectric excitation on the conservative-unforced system dynamics. The tags in the charts report the actuation voltages. (b) Comparison between the frequency response curve estimated with the present formulation (DPIM-P)
    and that predicted by simply projecting the piezoelectric force on the linear modal subspace obtained without taking into account the new static position (DPIM-A).}
    \label{fig:ccbeam_frf}
\end{figure}

Figure~\ref{fig:ccbeam_frf} collects the results of the simulations
performed considering the polarisation curves reported in Fig.~\ref{fig:polarisation_curves} a) for three different voltage bias
equal to $10$, $15$, and $20$~V, and put in evidence the expected strong hardening behaviour of the beam. 
The DPIM simulations using orders 7 and 6 for the autonomous and 
non-autonomous parts, respectively, are benchmarked
in Figure~\ref{fig:ccbeam_frf}a
against the full order HBFEM results with 7 harmonics,
showing a perfect agreement.
The parametrisation of the non-autonomous part has been performed  around the master mode eigenfrequency corresponding to the new fixed point i.e.\ $\Omega=5.399$\,rad/$\mu$s. 
An important remark is that the reduced model is obtained by considering only the lowest harmonic component of the forcing that resonates with the driven mode
(see  Eq.~\eqref{eq:decompforce}). 
It is anyway worth stressing that in specific applications 
higher-order harmonics of the forcing might induce parametric excitations that can be accounted for by the present formulation
\cite{thomas2013efficient,givois2020experimental,givois2021dynamics}. The backbone curves, also reported in Fig.~\ref{fig:ccbeam_frf}(a), show a marginal shift with increasing forcing amplitudes, as a consequence of the static deflection created by the non-zero mean value of the piezo actuation and illustrated in Fig.~\ref{fig:clamped_fixed}. 
Even though the frequency shift is tiny and almost negligible,
as expected from the small static displacements of Fig.~\ref{fig:clamped_fixed},  nevertheless it is taken into account in the procedure.

Concerning the numerical performance of the proposed formulation, we report that each FRC computed with the HBFEM takes around 12 hours on a standard workstation (Intel Xeon Gold 6140, 2.3 GHz, 128 GB RAM), 
while the DPIM approach takes 3 minutes to compute the parametrisation and few seconds to compute each frequency response curve using numerical continuation of periodic orbits with the MATCONT package \cite{dhooge2003matcont}.

Figure~\ref{fig:ccbeam_frf}(b) compares the outcomes of the complete DPIM procedure developed here (and hereafter labeled as DPIM-P for ``present formulation''), to a simplified one where two important assumptions routinely used in simulations have been considered. In this simplified version, denoted as DPIM-A for "approximate", the first assumption consists in taking into account the non-autonomous terms in a simplified manner, by simply projecting the piezo forces on the modal master modes. This assumption corresponds to a DPIM-$\mathcal{O}(p,0)$, {\em i.e.} an order 0 treatment of the external forcing, as commented for example in~\cite{opreni22high}, and used for instance in~\cite{touze2006nonlinear,vizzaccaro2021direct,vizza21high}. The second assumption used in DPIM-A consists in using the master eigenvectors of the unforced problem to perform the projection of the forcing, {\em i.e.} without taking into account the new static position and the shift of the eigenfrequency created by the constant forcing terms of the piezo. Saying things differently, one uses in DPIM-A the eigenvectors provided by the traditional stiffness matrix and not the tangent one $\tKu$, shown in Eq.~\eqref{eq:zitang00a}, since $\Uu_0$ is not considered. This assumption has been used in numerous examples in the past, see {\em e.g.}~\cite{frangi2020nonlinear,opreni2022fast,andrea112022}. For this specific example of the clamped-clamped beam with the selected polarisation, one can observe in Fig.~\ref{fig:ccbeam_frf}(b) a very slight difference between the two approaches, mainly due to the fact that the shift of the static position is negligible in the present case. This illustrates that simplified solutions can give, in many cases, correct results.

\begin{figure}[htb!]
    \centering
    \includegraphics[width = .99\linewidth]{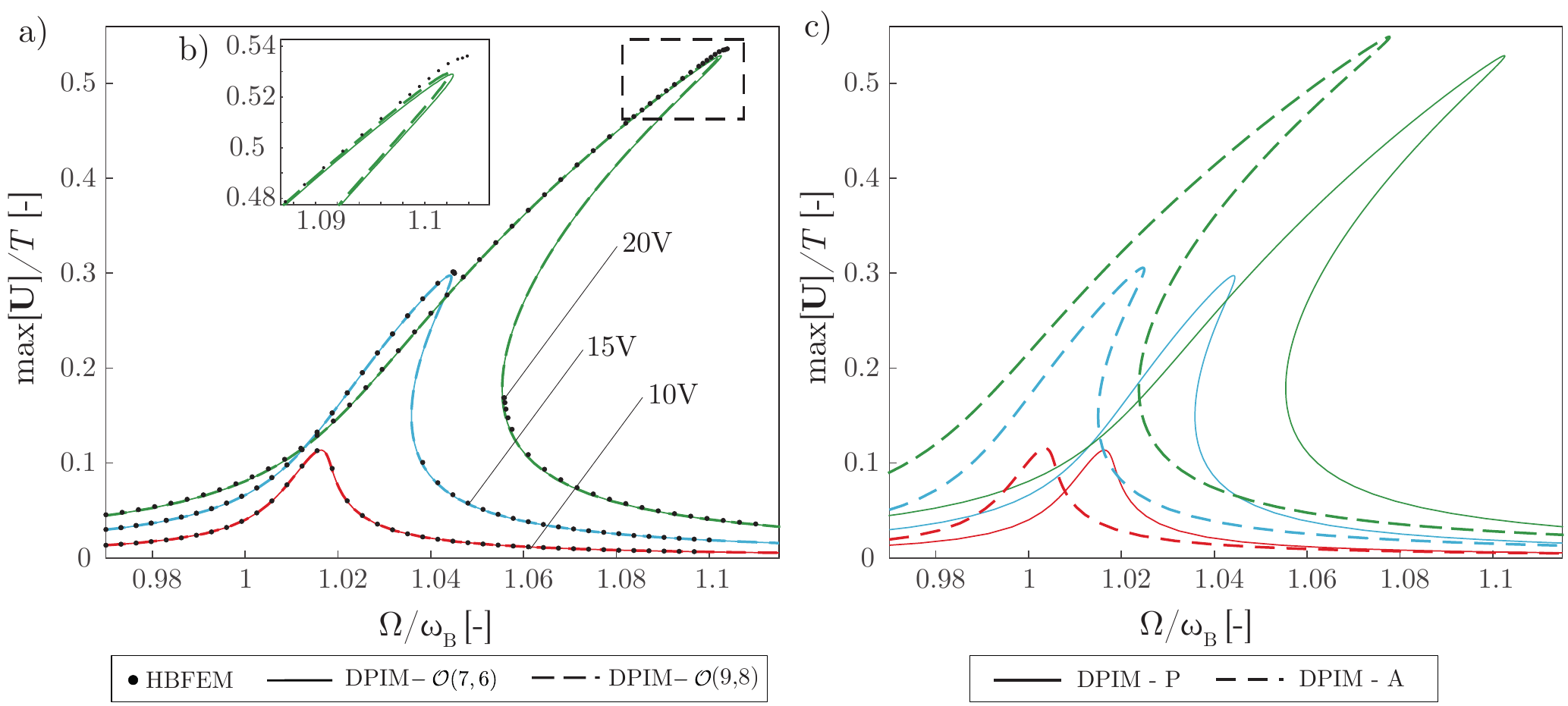}
    \caption{FRCs corresponding to different DPIM models.  (a) Comparison between the FRC estimated with the present formulation using two different approximation order $\mathcal{0}(7,6)$ and $\mathcal{0}(9,8)$.  (b) Enlarged view of the FRC peak.
    (c) Comparison between the FRC estimated with the presented formulation (DPIM-P) and the one predicted by simply projecting the piezoelectric force on the linear modal subspace with no shift of the fixed point  (DPIM-A).}
    \label{fig:FRF_order_cc_long}
\end{figure}

However, the difference between the two formulations becomes much more evident and important
when considering the modified polarisation curves reported in Fig.~\ref{fig:polarisation_curves_mod} to provide the actuation. 
These hysteresis loops are tailored to enhance the frequency shift introduced by the average piezoelectric forcing. 
The DPIM-P formulation yield an excellent match with HBFEM solutions plotted in Fig.\ref{fig:FRF_order_cc_long}(a), recovering the essential features of the FRCs up to large amplitudes: eigenfrequency shift, hardening behaviour and maximal amplitudes being fairly well reproduced.
Furthermore, the DPIM accuracy increases consistently with the expansion order, as underlined in the zoom Fig.\ref{fig:FRF_order_cc_long}(b), where one can observe the slight increase in accuracy obtained when moving from order $\mathcal{O}(7,6)$ to $\mathcal{O}(9,8)$. A small discrepancy at the very top of the FRC peak is nevertheless observed, pointing out that with this amplitude one starts to reach the accuracy limits provided by using a first-order approximation of the non-autonomous terms, see related discussions in~\cite{opreni22high}. 
The proposed model reproduces properly both the nonlinearity content and the piezoelectric-induced frequency shift as compared to the simplified formulation DPIM-A that only projects the piezoelectric force on the master mode and neglects the fixed point update, as illustrated in Fig.\ref{fig:FRF_order_cc_long}c). In that case, the change in the static position and in the linear eigenfrequency is too important, such that the two simplifying assumptions retained to build the model provided with DPIM-A do not hold anymore.


\begin{figure}[htb!]
    \centering
    \includegraphics[width = .99\linewidth]{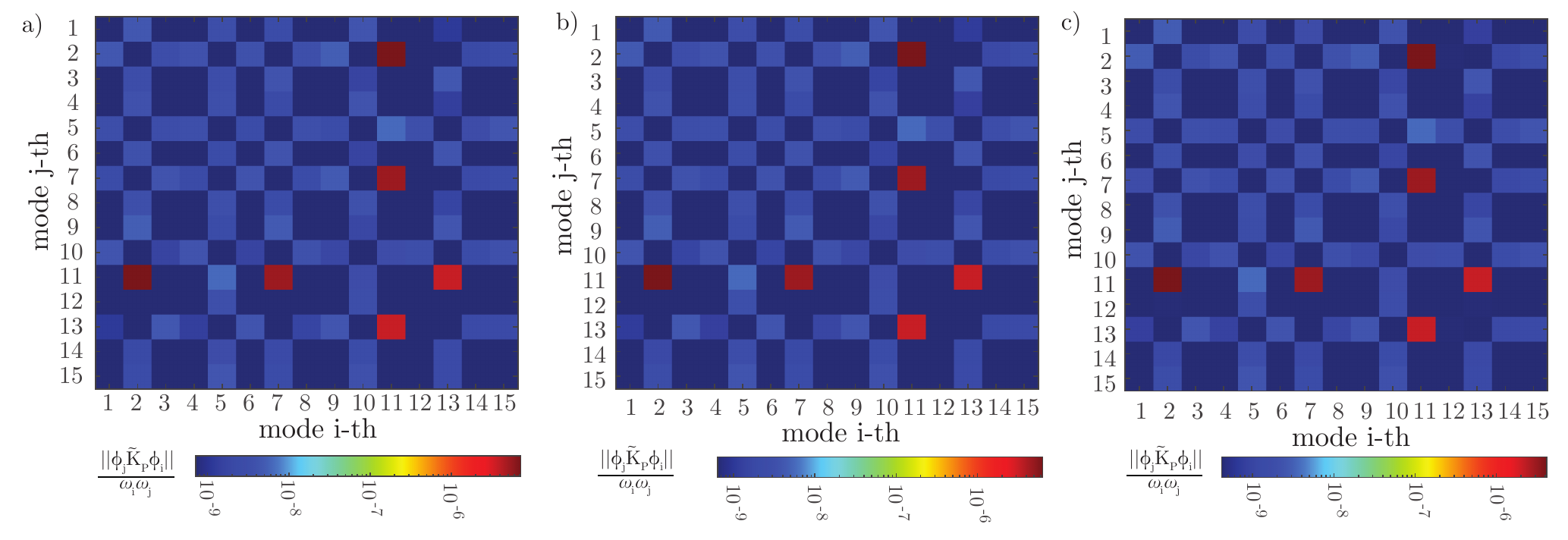}
    \caption{
    Visual representation of the influence of $\tKu_P$ on the first 15 eigenmodes. 
    Fig. a),b) and c) refer to the voltage levels 10V, 15V and 20V corresponding to the polarisation curves of Fig.~\ref{fig:polarisation_curves_mod}. Each coloured dot
    in the matrices
    represents  $|\bPhi^T_j\tKu_P\bPhi_i|/(\omega_i \omega_j)$}.
    \label{fig:kpiezo}
\end{figure}

To conclude the analysis of this academic example, the assumption that the added non-autonomous linear term on the right-hand side of Eq.~\eqref{eq:lin_mom_full_compact}, brings about negligible modifications to the eigenvalues, is verified. Indeed, the matrix $\tKu_P$ is expected to modulate in time the eigenvalues and eigenmodes of the system used to compute the parametrisation of master invariant manifolds. Consequently, this modulation needs to be small to avoid revising the whole computational scheme. Considering the first 15 eigenmodes, the $\tKu_P$ effect can be estimated by projecting the matrix onto the modal space. In particular, 
let us introduce the normalised projection of $\tKu_P$ onto the eigenvectors defined as $|\bPhi^T_j\tKu_P\bPhi_i|/(\omega_i \omega_j)$.
The resulting values are collected in the matrices depicted in Fig.~\ref{fig:kpiezo}. 
This representation highlights that the $\tKu_P$ introduces minor cross-couplings between some of the system modes as a consequence of the fact that the $\tKu_P$ matrix is not orthogonal to the eigenbases.
In particular, the master mode couplings, i.e. the first row and column in each figure, have a relative magnitude in the order of $10^{-8}$ at most.
This is reasonable and expected since the results are very accurate even with only one master mode. Furthermore, also the other cross-coupling terms are small with respect to the system stiffness, the largest terms having a relative magnitude of $10^{-5}$. 
This means that the change of non-autonomous eigenfrequencies is negligible, consistently with what remarked in \cite{opreni22high}.

The results on this academic example show a very good performance of the DPIM approach to offer predictive and accurate ROMs for a clamped-clamped beam actuated with different polarisation histories. Note that in Appendix~\ref{app:cantilever}, the case of a cantilever beam is analyzed. This benchmark is known to be more challenging for reduction methods mainly because of the folding of the invariant manifold corresponding to the fundamental bending mode at large amplitudes~\cite{vizza21high}. Anyway, results collected in Appendix~\ref{app:cantilever} show similarly an excellent behaviour of the DPIM.

\subsection{Micromirrors}
\label{sec:mirrors}

In the present Section we apply the reduction method
to MEMS micromirrors, a case of remarkable industrial interest
as they are key enabling components of many high-end
industrial applications.
The simulation of piezo micromirrors is a challenging task as even small nonlinearities in the Frequency Response Curves can degrade their optical performance significantly.
The analysis of their behaviour has recently stimulated intensive research \cite{opreni2021model,opreni22high,frangi2020nonlinear,opreni2021analysis}, highlighting the difficulty of generating accurate predictions in the presence of large rotations.
Actually, given the high dimensionality of finite element models of MEMS structures, 
the DPIM can be considered as the only available technique capable of exactly predicting the nonlinear dynamic response of MEMS components
within time spans compatible with the design and optimization of MEMS devices.

\begin{figure}[ht]
    \centering
    \includegraphics[width = .99\linewidth]{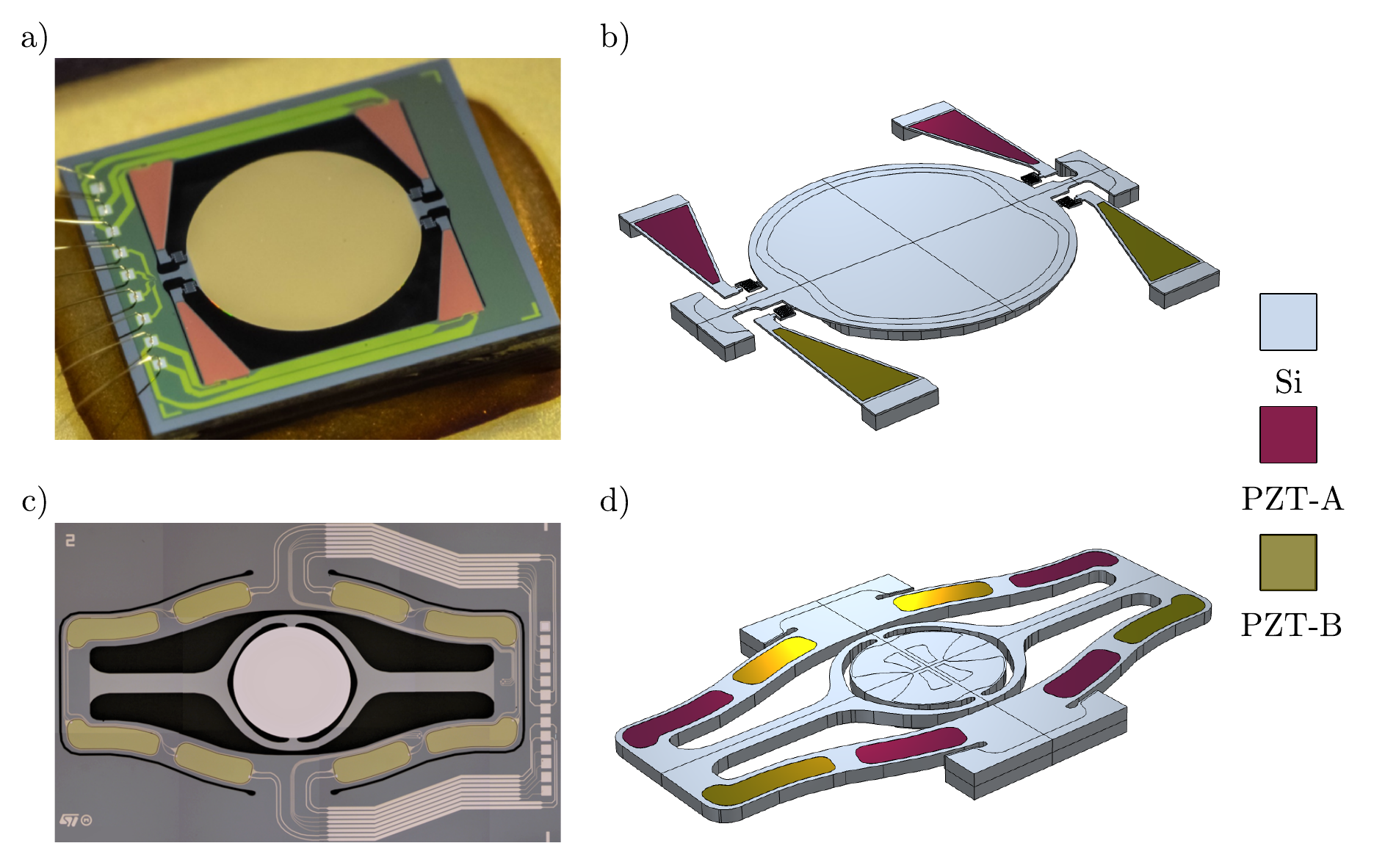}
    \caption{(a) Optical microscope image of Mirror A. (b) Geometry of the modeled device. Colors are used to distinguish the silicon structure of the device from the actuation patches, the latter being organised in groups A and B.}
    \label{fig:mirrors_img}
\end{figure}

The first mirror under consideration, hereafter labelled as Mirror A, is depicted in 
Fig.~\ref{fig:mirrors_img}(a) and has been fabricated through the P$\veps$Tra Thin-Film-Piezoelectric technology developed by STMicroelectronics. 
The yellow central circle is the reflective surface,
with a diameter of approximately 3000 $\mu$m and
a thickness of 20 $\mu$m, reinforced with a curvilinear
beam in order to minimize the dynamic deformation.
This surface is connected to the substrate via a pair of torsional springs.
The actuation is provided by PZT patches of thickness 2 $\mu$m,
visible in light orange, deposited on top of the four trapezoidal beams.
The polarisation histories for mirror A are plotted
in Figure~\ref{fig:polarisation_curves}a).
Due to the electrostrictive effect these bend activating the mirror rotation.
The actuation force is transmitted to the reflective surface through sets of folded springs.
As a first approximation the rotation axis of this mirror can be considered 
fixed, which allowed to derive a simple ROM in \cite{opreni2021analysis}
in terms of the rotation angle.

The micromirror is made of monocrystalline silicon with the [110]
orientation aligned with the torsional springs.
The materials properties are reported in Table \ref{tab:parameters_mir}.
The resonance frequency of this device is approximately 1950 Hz, up to imperfections in the fabrication process and the torsional mode has the lowest frequency 
in the spectrum.

\begin{table}[h]
	\begin{center}
		\begin{tabular}{|ccc|}
              \hline
              \multicolumn{3}{|c|}{PZT}\\
			\hline
			Coefficient & Value & Unit\\
			\hline
			$\cQ_{3333}$ 		& 0.097 	& m$^4$/C$^2$ 		\\
			$\cQ_{1133}$ 		& -0.046 	& m$^4$/C$^2$		\\
			$E$ 		& 70000	& MPa 		\\
			$\nu$ 		& 0.33 	& -		\\
			\hline
			\hline
                \multicolumn{3}{|c|}{Silicon}\\
                \hline                
			Coefficient & Value & Unit\\
			\hline
			$\cA_{1111}$ 		& 194250 	& MPa 		\\
   $\cA_{1122}$ 		& 35776 & MPa 		\\
   $\cA_{2222}$ 	 & 194250& MPa 		\\
   $\cA_{1133}$ 		& 64422 & MPa 		\\
   $\cA_{2233}$ 		& 64422 & MPa		\\
   $\cA_{3333}$ 		& 165605 & MPa 		\\
   $\cA_{2323}$ 		& 50591	& MPa 		\\
   $\cA_{3131}$ 		& 79237	& MPa 		\\
   $\cA_{1212}$ 		& 79237	& MPa 		\\
			\hline
		\end{tabular}
	\end{center}
	\caption{Constitutive parameters for PZT and silicon for the micromirrors examples. Single crystal silicon has an orthotropic mechanical behaviour, while PZT is considered isotropic.}
	\label{tab:parameters_mir}
\end{table}

The second micromirror, hereafter labelled Mirror B, features a different geometry. It is made of a reflective surface having a diameter of approximately 2000 $\mu$m
and a thickness of 150 $\mu$m.  This surface is connected to a gimbal structure as shown in Figure \ref{fig:mirrors_img}c. The gimbal is in turn anchored to ground. The resonance frequency of this device is approximately 25000 Hz. Actuation is obtained by means of eight PZT patches evidenced in Fig. \ref{fig:mirrors_img}d, organised in two groups and actuated with the same voltage laws given in Eq.~\eqref{eq:actuation_law}. 
The polarisation histories for mirror B are plotted
in Figure~\ref{fig:polarisation_curves}b).
Material properties and orientation are the same as for Mirror A (Table \ref{tab:parameters_mir}).
The simulation of this latter mirror poses important challenges,
as the torsional mode is only the fourth in the spectrum and is not well separated from other modes.
Globally, mirror B has a softening behaviour and it displays more clearly the classical unstable branch between the two saddle-node points
at any excitation voltage, whereas mirror A  displays hardening behaviour and a very short appearance of the unstable part for the largest amplitude tested. The axis of rotation is not clearly defined and 
the bulky reinforcement  induces a coupling with
translational motions, thus reducing the apparent stiffness.

\begin{figure}[ht]
    \centering
    \includegraphics[width = .99\linewidth]{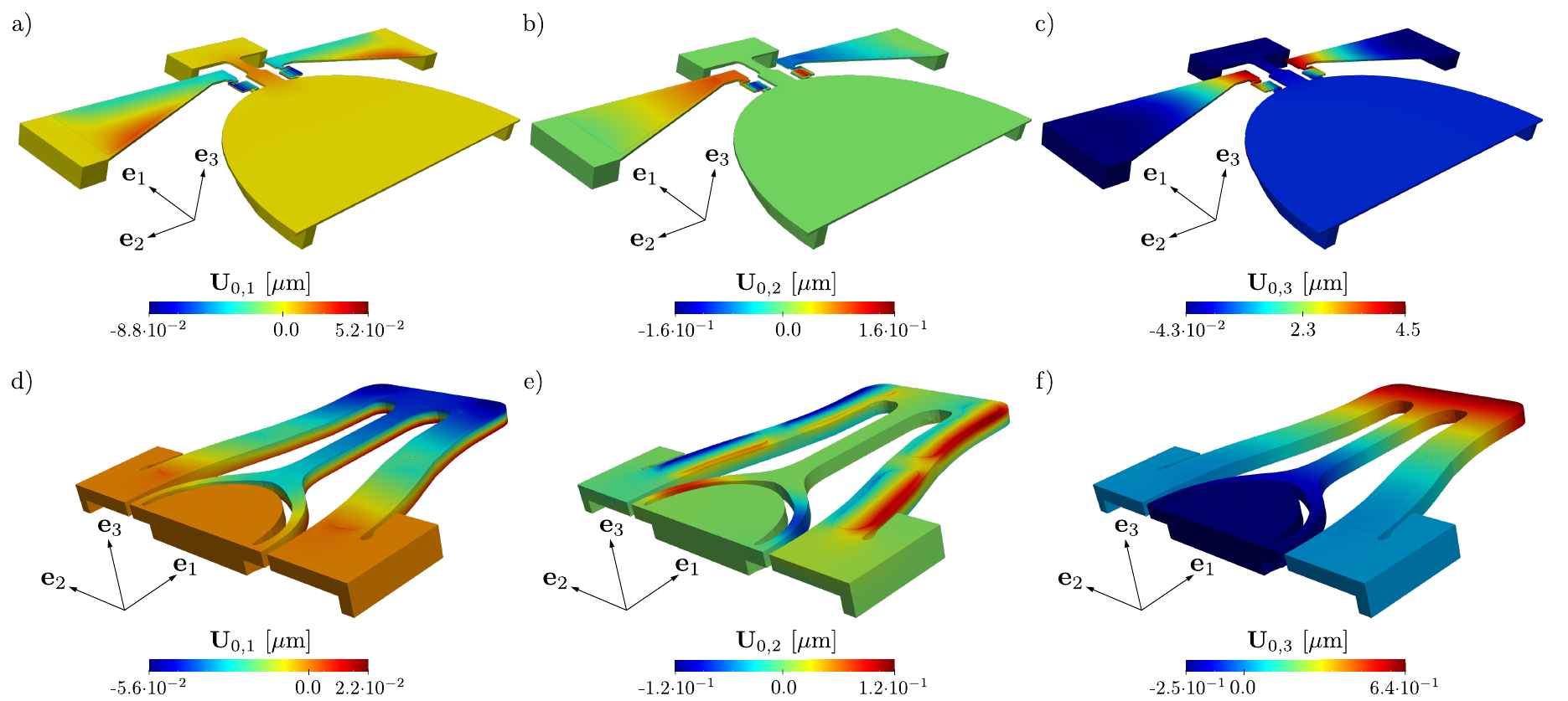}
    \caption{Displacement field components associated to the fixed point of the two micromirrors. In Figs. a), b) and c) are reported the static positions of Micromirror A for $V_0$ = 20 V  $\eu_1$,$\eu_2$,$\eu_3$ components respectively. In Figs. d), e) and f) are reported the static positions of Micromirror B for $V_0$ = 30 V  $\eu_1$,$\eu_2$,$\eu_3$ components respectively.}
    \label{fig:static_mirror}
\end{figure}

As highlighted in the previous example, the applied piezoelectric voltage bias induces initial stresses in the structure that alter the fixed point and the corresponding eigenfrequencies of the system. The static displacement fields associated with the maximum bias for the two micromirrors are reported in 
Fig.~\ref{fig:static_mirror}. Since inelastic strains are acting on portions of the structures that can freely displace, no relevant changes of the torsional mode of the devices are expected, as also highlighted in the upcoming results.

Preliminary analyses in \cite{opreni2021model,vizza21high}
have shown that a high order DPIM expansion is required in this specific application.
The first validation of the method is performed by comparing the results of the reduced model with the HBFEM simulations implementing the same formulation. 
Results are reported for voltage amplitudes equal to $10$, $15$, and $20$ V for Mirror A, and $20$, $25$, and $30$ V for Mirror B. These values are comparable with those applied during the real functioning of the device.

While the devices analysed in \cite{opreni2021analysis} are exactly those presented in this work, a coarser finite element discretisation is 
here adopted to reduce the computational burden for full order simulations. The mesh for Mirror A consists of 15341 nodes and the Fourier expansion used to approximate the sinusoidal motion of the device is taken of order 5, i.e.\ the resulting number of degrees of freedom in Fourier domain is equal to 506253. For Mirror B, a discretisation based on 21260 nodes is used with a Fourier expansion of order 7 that results in 956700 degrees of freedom in the Fourier domain. 

\begin{figure}[ht]
    \centering
    \includegraphics[width = .99\linewidth]{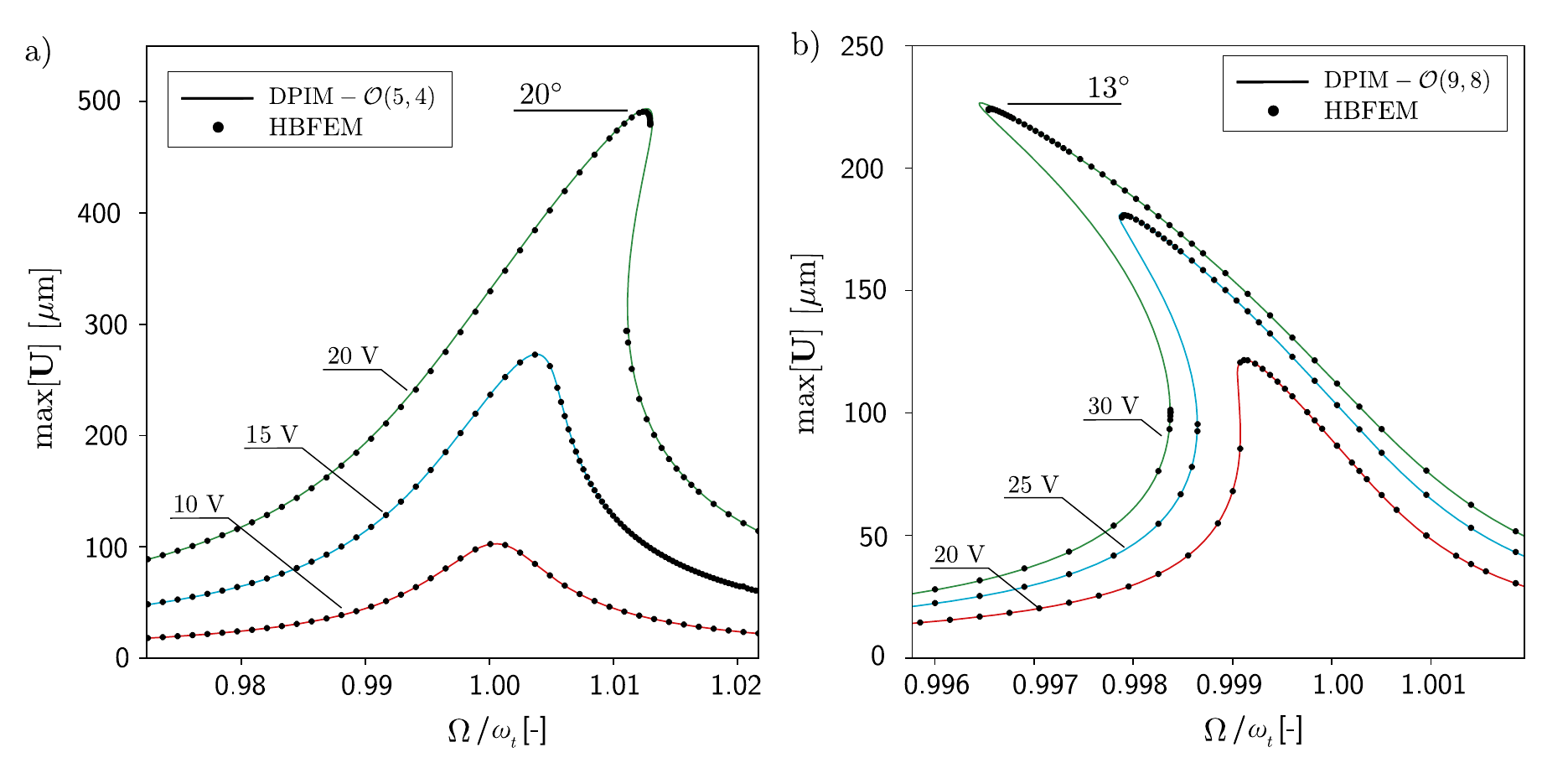}
    \caption{Comparison between numerical results computed from the direct parametrisation method for invariant manifolds and full order HBFEM simulations. (a) reports the results computed for Mirror A, while (b) reports the results computed for Mirror B. In both charts the corresponding voltage amplitudes $V_0$ are reported as tags, and the same is done for the maximum rotation amplitude reached by the two devices.}
    \label{fig:lynx_perseus_results}
\end{figure}

A comparison between full order HBFEM simulations and the reduced model obtained from the direct parametrisation method for invariant manifolds is reported for both devices in Fig. \ref{fig:lynx_perseus_results}, where Rayleigh damping parameters are set as $\alpha=\omega_{t}/100$  and $\beta=0$ for Mirror A and $\alpha=\omega_{t}/1000$ and $\beta=0$ for Mirror B, with $\omega_{t}$ frequency of the torsional mode. 
These $\alpha$ values yield quality factors compatible with 
those extracted from experimental data using the relations provided by Davis for cubic oscillators \cite{davis2011measuring}. We remark that damping models for MEMS systems operating at high amplitudes are often nonlinear as a result of convective effects \cite{di2022arbitrary}. However, parameter identification for nonlinear models is often impractical and often unnecessary given uncertainties in the measurement setups, so a linear model is often preferred.

\begin{figure}[ht]
    \centering
    \includegraphics[width = .99\linewidth]{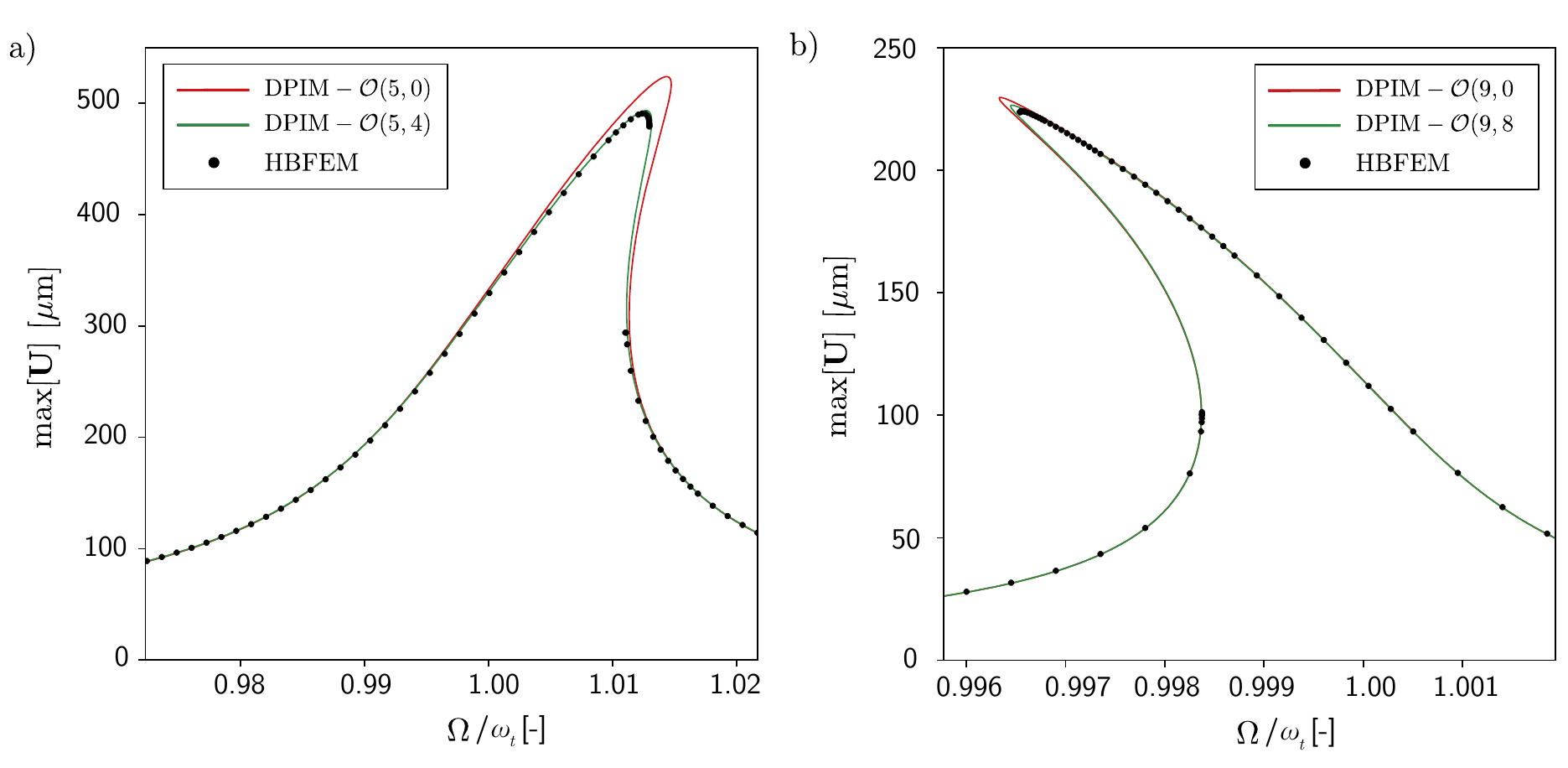}
    \caption{Convergence of the non-autonomous part of the DPIM when increasing the expansion order of the $\veps^1$-invariance equation for Mirror A (a) and Mirror B (b).}
    \label{fig:lynx_perseus_results_order}
\end{figure}

The charts highlight a perfect agreement between the two numerical models both in terms of nonlinearity of the curves and maximum amplitude. 
A further remark on Fig. \ref{fig:lynx_perseus_results} is the difference in expansion order adopted for the two devices. Mirror A achieves convergence already at order 5 of the asymptotic expansion. This has already been evidenced in \cite{opreni2021model}, where a low order direct normal form approach was used and excellent results in terms of nonlinearity have been obtained. This is consistent with the structure geometry. Indeed, Mirror A features a small reinforcement beneath its reflective surface that does not alter the kinematics of the structure. As a result, its behaviour is similar to that of a flat structure as a cantilever, with a mild hardening response. This implies that a low order parametrisation procedure that accounts for non-resonant stretching modes is sufficient to correctly predict the nonlinear dynamic response of the structure. On the other hand, Mirror B features a bulky support beneath the reflective surface. This introduces eccentricity in the structure and as a result strong quadratic coupling between the out-of-plane mode and the torsional mode of the structure. The consequence is that the order of expansion of the method to achieve the correct nonlinearity trend increases. As already evidenced in \cite{vizza21high}, for the same device geometry, an order 7 asymptotic expansion is necessary to correctly catch the nonlinear response of the device. Morevover, in the present work an order 9 expansion of the $\veps^0$-invariance equation is proposed and the correct nonlinearity type is predicted as highlighted in Fig. \ref{fig:lynx_perseus_results}(b). 

A further check is the convergence of the method with the expansion order of the $\veps^1$ terms. Results are reported in Fig. \ref{fig:lynx_perseus_results_order}, where comparisons between reduced models obtained with either high order or zero order expansion of the non-autonomous invariance, are given. Results highlight the necessity to adopt a high order parametrisation order for both autonomous and non-autonomous invariance equations to achieve the same results predicted by full order HBFEM simulations. In particular, we highlight that a high order expansion of the $\veps^1$-invariance is necessary in presence of strong changes in configuration. Indeed, as it can be appreciated in Fig. \ref{fig:lynx_perseus_results}, Mirror A is subjected to rotations up to 20$^{\circ}$. On the other hand, Mirror B is subjected to rotations only up to 12$^{\circ}$ hence, as shown in  
Fig.~\ref{fig:lynx_perseus_results_order}, the advantage of adopting a high order parametrisation for the non-autonomous part is less evident. 

\begin{figure}[htb!]
    \centering
    \includegraphics[width = .99\linewidth]{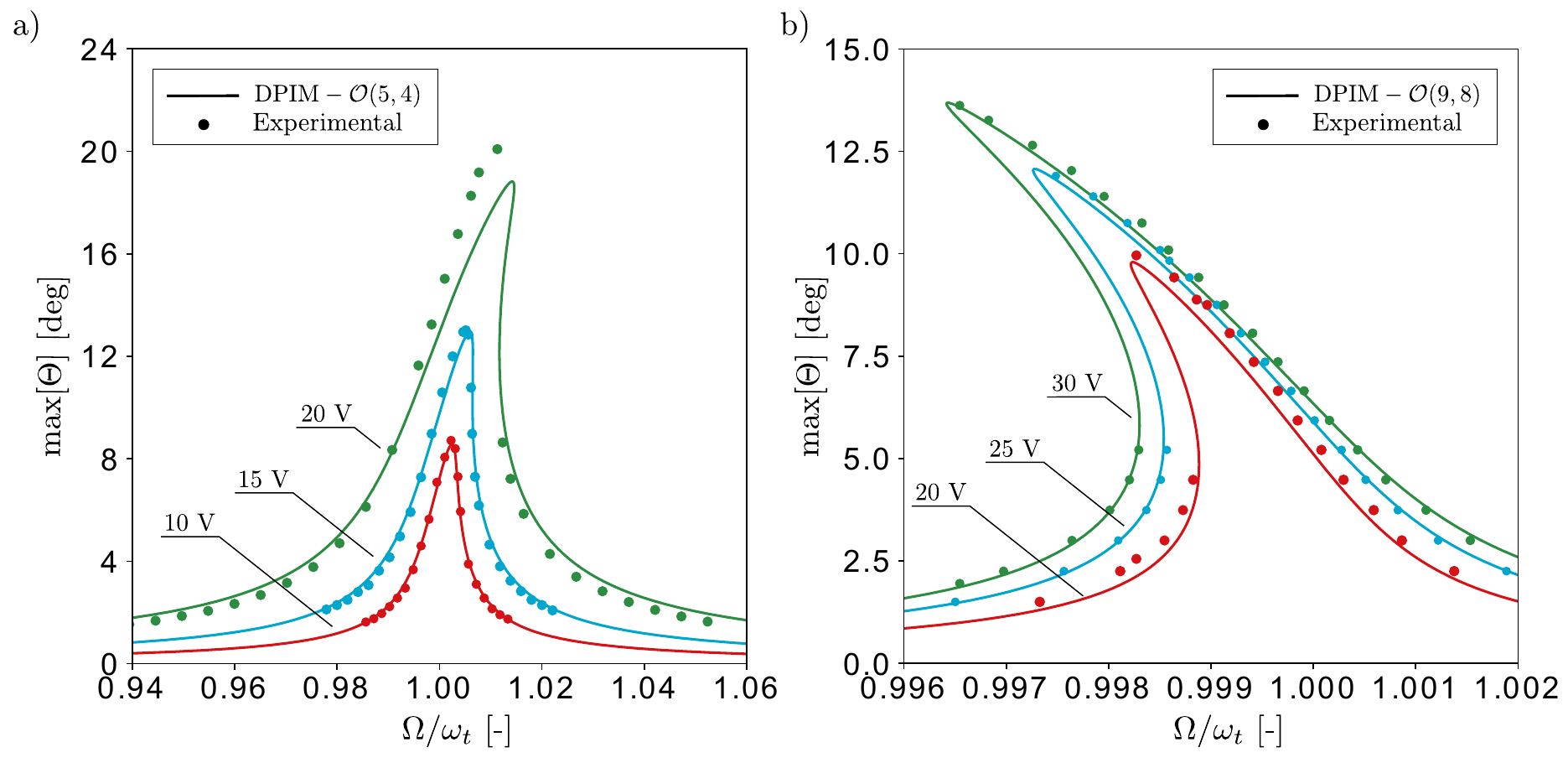}
    \caption{Comparison between experimental and numerical FRCs for Mirror A (a) and Mirror B (b).}
    \label{fig:lynx_perseus_experiments}
\end{figure}

Finally, an experimental validation of the predictions given by the model is reported
in Fig. \ref{fig:lynx_perseus_experiments}. 
Experimental data are identical to those detailed in \cite{opreni2021analysis}, where the setup and device characterisation procedures are discussed. Experiments have been conducted in frequency control, hence no information regarding unstable branches is available.
Overall, the accuracy of the method is good in terms of both nonlinearity of the predicted response and absolute value of the oscillation amplitude, hence highlighting the valuable predictive capabilities of the method. The discrepancies observed on mirror A at high voltages and large apertures can be ascribed to the complex flow patterns that develop around the mirror at large amplitude thus inducing strongly nonlinear damping.
Here, on the contrary, a linear Rayleigh damping model has been assumed with a constant value of the $\alpha$ coefficient for all the analyses run on a given mirror,
coherently with the results reported in Figures~\ref{fig:lynx_perseus_results}.
It is worth mentioning that in \cite{opreni2021analysis} a better agreement 
has been obtained by using different values of the quality factor
for each voltage value, calibrated from experiments and ad-hoc formulas \cite{davis2011measuring}. This could have been reproduced here to show the best possible fitting, nevertheless the focus has been set on showing the numerical predictions given by the ROM for a single set of parameters that are not fitted when changing the excitation frequency or amplitude.

Concerning the numerical performance of the proposed formulation,  the DPIM approach 
takes 2 minutes for micromirror 1, and 8 minutes for micromirror 2 to compute the whole parametrisation, on the other hand an HBFEM simulation requires at least 12 hours to compute only 20 pints of the FRFs on a workstation with Intel Xeon Gold 6140, 2.3 GHz, 128 GB RAM.

\section{Conclusions}

The Direct Parametrisation for Invariant Manifolds (DPIM)
approach has been extended to piezo micro actuators, accounting for the specificities introduced by the piezo forcing. 
This investigation thus brings this recently developed powerful technique 
to a new level of maturity and practical utility.

Assuming that the polarisation history is known from experimental measurements,
the proposed approach accounts exactly for geometrical, inertia and material nonlinearities induced by the Landau-Devonshire constitutive modelling 
of electrostrictive effects in ferroelectrics.
The inclusion of time dependent and configuration dependent piezo forcing
has required to modify the DPIM procedures recently published,
which were limited to deformation independent forcings.
Both the treatment of the autonomous and non-autonomous parts have been detailed
exploiting the similarities with previous developments and permitting efficient coding.  
It has also been shown that an accurate ROM can be obtained by deriving the non-autonomous terms for a single, 
fixed external frequency, thus importantly reducing the computational burden.

This simulation tool developed represents a major achievement, as 
it provides ideal simulation capabilities for a whole class of applications involving
piezo actuators.
Starting from the full 3D FEM model of the device, we have shown how to compute at the same time both the non-linear mappings for displacement and velocity nodal values and the reduced dynamics. 
This is a major benefit of the proposed approach as
modern MEMS often cannot be efficiently described resorting to simplified structural theories.

Even if the focus has been set on  piezo-MEMS actuators fabricated using Lead Zirconate Titanate (PZT), which is deposited in the form of a thin film sol-gel on bulk silicon, the proposed approach can be easily extended to other piezo materials.
The tool has been benchmarked on academic and industrial applications, both against
an accurate full-order Harmonic Balance Method and experiments, 
thus validating the whole set of underlying assumptions.

Even if we have not addressed here a fully coupled multi-physics problem,
this investigation represents a first step in this direction and opens  
a fascinating challenge in the field of reduced order modelling.

\section*{Data availability}
The data that support the findings of this study are available on request from the corresponding author.

\section*{Declaration of Competing Interest}
The authors declare that they have no known competing financial interests or personal relationships that could have appeared to influence the work reported in this paper.

\bibliographystyle{unsrt}
\bibliography{biblio}

\appendix


\section{Finite element discretisation of the governing equations}
\label{app:FEM_formulation}

In the present work Eq. \eqref{eq:lin_mom_weak_full_explicit} is numerically discretised using a Bubnov-Galerkin finite element formulation. To thos aim, let us introduce a finite element discretisation of the type:
\begin{equation}
    \uu \approx \uu_h = \sum_{i=1}^{N}\lu_i U_i, \qquad \wu \approx \wu_h = \sum_{i=1}^{N}\lu_i W_i,
\end{equation}
with $\lu_i$ nodal basis functions and $U_i$, $W_i$ nodal values of displacement field and test function, respectively. Upon substitution of the finite element discretised field in Eq. \eqref{eq:lin_mom_weak_full_explicit}, we can derive the following discretised quantities:

\begin{subequations}
  \begin{align}
    & \int_{B} \rho\,\ddot{\uu}_h \cdot\wu_h\,\dd B = \sum_{i,j=1}^{N} W_i \int_{B} \rho\,{\lu}_j \cdot\lu_i\,\dd B\, \ddot{U}_j = \Wu^{T} \Mu \ddot{\Uu}, \\
    & \int_{B} \sym(\nabla\uu_h):\cA:\sym(\nabla\wu_h)\,\dd B = \sum_{i,j=1}^{N} W_i \int_{B} \sym(\nabla\lu_j):\cA:\sym(\nabla\lu_i)\,\dd B\, U_j = \Wu^{T} \Ku \Uu,\\
    & \int_{B} \sym(\nabla\uu_h):\cA:\sym(\nabla^{T}\wu_h\cdot \nabla\uu_h)\,\dd B + \frac{1}{2} \int_{B} \sym(\nabla\wu_h):\cA:\sym(\nabla^{T}\uu_h\cdot \nabla\uu_h)\,\dd B = \nonumber \\
    & \sum_{i,j,k=1}^{N} W_i  \int_{B} \sym(\nabla\lu_j):\cA:\sym(\nabla^{T}\lu_i\cdot \nabla\lu_k)\,\dd B + \frac{1}{2} \int_{B} \sym(\nabla\lu_i):\cA:\sym(\nabla^{T}\lu_j\cdot \nabla\lu_k)\,\dd B \, U_j U_k = \Wu^{T}\G(\Uu,\Uu),\\
    & \frac{1}{2} \int_{B} \sym(\nabla^{T}\uu_h\cdot \nabla\uu_h):\cA:\sym(\nabla^{T}\wu_h\cdot \nabla\uu_h)\,\dd B = \nonumber  \\
    & \sum_{i,j,k,l=1}^{N} \frac{1}{2} W_i \int_{B} \sym(\nabla^{T}\lu_j\cdot \nabla\lu_k):\cA:\sym(\nabla^{T}\lu_i\cdot \nabla\lu_l)\,\dd B \, U_j U_k U_l = \Wu^{T}\H(\Uu,\Uu,\Uu), \\
    & \int_{B} \Su^p:\sym(\nabla\wu_h) \,\dd B = \sum_{i=1}^{N} W_i \int_{B} \Su^p:\sym(\nabla\lu_i) \,\dd B  = \Wu^{T} \Fu_P, \\
    & \int_{B} \Su^p:\sym(\nabla^{T}\uu_h\cdot\nabla\wu_h) \,\dd B = \sum_{i,j=1}^{N} W_i \int_{B} \Su^p:\sym(\nabla^{T}\lu_j \cdot\nabla\lu_i) \,\dd B\, U_j = \Wu^{T} \Ku_P \Uu,
  \end{align}
\end{subequations}

where the last two expressions represent the additional terms that stem from the employed piezoelectric formulation. We remark that the same expression holds in presence of a generic pre-stress term. Furthermore, the applicability of the presented method is not affected by the choice of the numerical scheme adopted to discretise the governing equations.

\section{Solution scheme for the static analysis}
\label{app:static_analysis}

Solution of Eq. \eqref{eq:static} is performed using a standard Newton-Raphson procedure. To this aim, let us report the partial differential equation that yields such system:

\begin{equation}\label{eq:ppv_static}
    \int_{B} \eu[\uu,\uu]:\cA:\delta\eu[\uu,\wu] \, \dd B = \int_{B} \Su^p_0:\delta\eu[\uu,\wu]\, \dd B, \forall \,\wu \in\cC(\zerou)
\end{equation}

with $\Su^p_0$ mean value of the inelastic stress caused by electrostriction. Square brackets are used specify on which variables each operator acts. Linearisation of Eq. \eqref{eq:ppv_static} around of generic known configuration $\bar{\uu}$ yields the following:

\begin{align}\label{eq:ppv_linearised}
    & \int_{B} \eu[\bar{\uu},\bar{\uu}]:\cA:\delta\eu[\bar{\uu},\wu] \, \dd B + \int_{B} \delta\eu[\bar{\uu},\delta\uu]:\cA:\delta\eu[\bar{\uu},\wu] \, \dd B + \int_{B} \eu[\bar{\uu},\bar{\uu}]:\cA:\delta\eu[\delta\uu,\wu] \, \dd B = \nonumber \\
    & \int_{B} \Su^p_0:\delta\eu[\bar{\uu},\wu] \, \dd B + \int_{B} \Su^p_0:\delta\eu[\delta\uu,\wu] \, \dd B, \quad \forall\, \wu \in \cC(\zerou),
\end{align}

where $\delta\uu$ is the displacement variation with respect to $\bar{\uu}$. We highlight that Eq. \eqref{eq:ppv_linearised} is linear with respect to the displacement increment $\delta\uu$. As a result, we can collect term that depend linearly on the displacement $\delta\uu$ at the left hand side and terms that can be computed from the configuration $\bar{\uu}$ on the right hand side:

\begin{align}\label{eq:ppv_linearised_collect}
    & \int_{B} \delta\eu[\bar{\uu},\delta\uu]:\cA:\delta\eu[\bar{\uu},\wu] \, \dd B + \int_{B} \eu[\bar{\uu},\bar{\uu}]:\cA:\delta\eu[\delta\uu,\wu] \, \dd B - \int_{B} \Su^p_0:\delta\eu[\delta\uu,\wu] \, \dd B =  \nonumber \\
    & - \int_{B} \eu[\bar{\uu},\bar{\uu}]:\cA:\delta\eu[\bar{\uu},\wu] \, \dd B +
     \int_{B} \Su^p_0:\delta\eu[\bar{\uu},\wu] \, \dd B , \quad \forall\, \wu \in \cC(\zerou),
\end{align}

where the integrals on the left hand side are in order the material tangent, the geometrical tangent, and the piezoelectric force tangent. The right hand side is simply the residual of the equation. Upon finite element discretisation of Eq. \eqref{eq:ppv_linearised_collect} it is possible to solve it iteratively by updating the configuration $\bar{\uu}$ until it minimises the total energy defined for Eq.  \eqref{eq:ppv_static}.

\section{Detailed explicit expressions of the ROM}
\label{app:ROMdetail}

This appendix is devoted to emphasizing some calculation details related to the derivation of the reduced-order models using the direct parametrisation of invariant manifolds (DPIM). The main developments have been reported in~\cite{vizza21high,opreni22high} for autonomous and non-autonomous systems encompassing geometric nonlinearity. Here the general developments are adapted to tackle the present problem with the piezoelectric forcing terms. A special emphasis is put on the additional treatments needed to take into account the new terms.

One of the main feature of the method is to derive arbitrary order homological equations, expressed at the specific level of a given monomial. To that purpose, all unknown functions (nonlinear mappings and reduced dynamics) are expanded, and each order-$p$ term is developped according to the summation of monomials written as:

\begin{subequations}\label{eq:polyexpandmon_HO}
\begin{align}
\P{\WUfun} = &\, 
\sum_{i_1=1}^{2n} \sum_{i_2=1}^{2n} \ldots \sum_{i_p=1}^{2n} 
\p{\WU_{\lbrace i_1 i_2 \ldots i_p\rbrace}} \; \zr_{i_1} \zr_{i_2} \ldots \zr_{i_p} = 
\sum_\Is \p{\WU_\Is} \; \p{\pi_\Is},
\\
\P{\WVfun} = &\, 
\sum_{i_1=1}^{2n} \sum_{i_2=1}^{2n} \ldots \sum_{i_p=1}^{2n} 
\p{\WV_{\lbrace i_1 i_2 \ldots i_p\rbrace}} \; \zr_{i_1} \zr_{i_2} \ldots \zr_{i_p} = 
\sum_\Is \p{\WV_\Is} \; \p{\pi_\Is},\\
 \P{\ffun} = &\,
\sum_{i_1=1}^{2n} \sum_{i_2=1}^{2n} \ldots \sum_{i_p=1}^{2n} 
\p{\fu_{\lbrace i_1 i_2 \ldots i_p\rbrace}} \; \zr_{i_1} \zr_{i_2} \ldots \zr_{i_p}
= \sum_\Is \p{\fu_\Is} \; \p{\pi_\Is}.
\end{align}
\end{subequations}

In these equations, following~\cite{vizza21high,opreni22high}, a given monomial is represented with $\p{\pi_\Is}$ expressed as:

\begin{equation}
    \p{\pi_\Is} = z_{i_1}z_{i_2}\,...z_{i_p},
\end{equation}

which is an order-$p$ monomial of the normal coordinates. For the ease of treatments, each monomial can be simply tracked by the set of indices $\Is$ as:

\begin{equation}
    \Is = \{i_1 i_2 \,... i_p\}.
\end{equation}

The cardinal number of $\Is$  equals p, meaning that in this ordering, indices with multiplicity higher than one are  repeated.

These solution forms for the unknowns mappings and reduced dynamics are substituted into Eqs. \eqref{eq:invariance_zero} and \eqref{eq:invariance_one_proj}, that are respectively the invariance equation at order $\varepsilon^0$ for the autonomous problem, and invariance equation at order $\varepsilon^1$ for the non-autonomous terms. Upon collecting and identifying same powers for each index set $\Is$ , then one is able to write  the homological equations that need to be solved to retrieve mappings and reduced dynamics of the reduced model. For the autonomous case, Eq.~\eqref{eq:invariance_zero}, one is thus able to write the homological equation at any order $p$, and for an arbitrary monomial $\Is \in \cH^{(p)}$, where $\cH^{(p)}$ refers to the set of all combination of indices with degree $p$:
\begin{subequations}
\begin{align}
	& \forall\, p = 2,...,o, \quad \forall\, \Is \in \cH^{(p)},  \nonumber\\
&
\Mu \p{\WV_\Is}\sigma_\Is 
+ \sum_{s=1}^{2n}\left(\Mu\phiu_s \lambda_s \,\p{\fr_{s\,\Is}} \right)
+ \Mu \p{\FV_\Is}
+ \Cu \p{\WV_\Is}
+ \Ku_{T} \p{\WU_\Is} + \p{\FG_\Is} + \p{\FH_\Is} = \zerou,
\\
& 
\Mu \p{\WU_\Is}\sigma_\Is 
+ \sum_{s=1}^{2n}\left(\Mu\phiu_s  \p{\fr_{s\,\Is}}\right)
+ \Mu\p{\FU_\Is}-\Mu\p{\WV_\Is}\p{\pi_\Is}
=
\zerou,
\label{eq:homo_eqs_HO}
\end{align}
\end{subequations}

where the new introduced quantities are defined as:

\begin{subequations}\label{eq:veps_zero_homological}
  \begin{align}
    \sigma_{\Is} =&\, \lambda_{i_1} + \lambda_{i_2} + ... + \lambda_{i_p}, \label{eq:veps_zero_homological-a}\\
\p{\FU_\Is} =&\,  
\sum_{s=1}^{2n}\sum_{k=2}^{p-1}\sum_{l=0}^{p-k}
\WU^{(p-k+1)}_{\lbrace i_1 \ldots i_l s\, i_{l+k+1} \ldots i_p \rbrace}
\fr^{(k)}_{s\lbrace i_{l+1} \ldots i_{l+k} \rbrace},
\label{eq:veps_zero_homological-b_NA}\\ 
\p{\FV_\Is} =&\,  
\sum_{s=1}^{2n}\sum_{k=2}^{p-1}\sum_{l=0}^{p-k}
\WV^{(p-k+1)}_{\lbrace i_1 \ldots i_l s\, i_{l+k+1} \ldots i_p \rbrace}
\fr^{(k)}_{s\lbrace i_{l+1} \ldots i_{l+k} \rbrace}. \label{eq:veps_zero_homological-c_NA}\\
\p{\FG_\Is} =&\, 
\sum_{k=1}^{p-1}
\tG(\WU^{(k)}_{\lbrace i_1 \ldots i_k\rbrace},
\WU^{(p-k)}_{\lbrace i_{k+1} \ldots i_p\rbrace}),
\label{eq:g_explicit_HO}\\
\p{\FH_\Is} =&\, 
\sum_{k=1}^{p-2}\sum_{l=1}^{p-k-1}
\H(\WU^{(k)}_{\lbrace i_1 \ldots i_k \rbrace},
\WU^{(l)}_{\lbrace i_{k+1} \ldots i_{k+l} \rbrace},
\WU^{(p-k-l)}_{\lbrace i_{k+l+1} \ldots i_p \rbrace}) .\label{eq:h_explicit_HO}
  \end{align}
\end{subequations}

A similar treatment is operated for the non-autonomous problem at order $\varepsilon^1$, and again the homological equations can be written for arbitrary order $p$ and for any monomial designated with the set of indices $\Is$ as:

\begin{subequations}\label{eq:zihomologicpNonAUt_NA}
  \begin{align}
	& \forall\, p = 0,...,q, \quad \forall\, \Is \in \cH^{(p)},  \nonumber\\
	& \hsigma_{\Is}\Mu \hWV^{(p)}_{\Is} + \Cu \hWV^{(p)}_{\Is} + \Ku_{T} \hWU^{(p)}_{\Is} + \sum_{s=1}^{2n} \left(\lambda_{s} \hf_{s\,\Is}^{(p)} \Mu\phiu_{s} \right) + \Mu \hFV_{\Is}^{(p)}  + \hFG^{(p)}_{\Is} + \hFH^{(p)}_{\Is}  = \hFu^{(p)}_{\Is}, \\
	& \hsigma_{\Is}\Mu \hWU^{(p)}_{\Is} - \Mu \hWV_{\Is}^{(p)} + \sum_{s=1}^{2n} \left(\hf_{s\,\Is}^{(p)} \Mu\phiu_{s}  \right) + \Mu \hFU_{\Is}^{(p)}  = \zerou,
  \end{align}
\end{subequations}

where all quantities not previously defined are given as:

\begin{subequations}\label{eq:zinonautgene_NA}
  \begin{align}
	\hsigma_{\Is} =&\, \hlm + \lambda_{i_1} + \lambda_{i_2} + ... + \lambda_{i_p}, \label{eq:zinonautgenea_NA}\\
	\hFU_{\Is}^{(p)} = &\, \sum_{s=1}^{2n}  \left( \sum_{k=2}^{p} \sum_{l=0}^{p-k} \hWU^{(p-k+1)}_{\{i_1 ... i_{l}\,s\,i_{l+k+1} ... i_{p}\}} f_{s\{i_{l+1}...i_{l+k}\}}^{(k)} + \right. \nonumber \\
	&\, \bigg. \sum_{k=0}^{p-1} \sum_{l=0}^{p-k} \WU^{(p-k+1)}_{\{i_1 ... i_{l}\,s\,i_{l+k+1} ... i_{p}\}} \hf_{s\{i_{l+1}...i_{l+k}\}}^{(k)} \bigg),\\
    \hFV_{\Is}^{(p)} = &\, \sum_{s=1}^{2n}  \left( \sum_{k=2}^{p} \sum_{l=0}^{p-k} \hWV^{(p-k+1)}_{\{i_1 ... i_{l}\,s\,i_{l+k+1} ... i_{p}\}} f_{s\{i_{l+1}...i_{l+k}\}}^{(k)} + \right. \nonumber \\
	&\, \bigg. \sum_{k=0}^{p-1} \sum_{l=0}^{p-k}  \WV^{(p-k+1)}_{\{i_1 ... i_{l}\,s\,i_{l+k+1} ... i_{p}\}} \hf_{s\{i_{l+1}...i_{l+k}\}}^{(k)} \bigg),\\
	\hFG_{\Is}^{(p)} =&\, 2 \sum_{k=1}^{p} \tG\left( \WU^{(k)}_{\{i_1...i_k\}}, \hWU^{(p-k)}_{\{i_{k+1}...i_p\}}   \right), \label{eq:GT_na}\\
	\hFH_{\Is}^{(p)} =&\, 3 \sum_{k=1}^{p-1}\sum_{l=1}^{p-k} \H\left( \WU^{(k)}_{\{i_1...i_k\}}, \WU^{(l)}_{\{i_{k+1}...i_{k+l}\}} , \hWU^{(p-k-l)}_{\{i_{k+l+1}...i_{p}\}}   \right),\\
	\hFu^{(0)}_{\Is} =&\, \tilde{\Fu}_{P}, \quad \mbox{and} \quad \hFu^{(p)}_{\Is} =\tilde{\Ku}_{P}\WU^{(p)}_{\{i_1\,...\,i_p\}}, \quad \forall \; p \geq 1. \label{eq:zinonautgenef_NA}
  \end{align}
\end{subequations}

It is possible to notice that the structure of the resulting homological equations is identical to that already derived in \cite{vizza21high,opreni22high}, the only differences being brought by the definition of the stiffness matrix, the quadratic operator, and the right-hand side of the non-autonomous homological equations. Indeed, the stiffness operator is given by the tangent stiffness matrix computed in one of the fixed points of the system. The shift of the fixed point with respect to the origin changes the quadratic operator, as highlighted in Eq. \eqref{eq:GT_na}. Finally, we report that the last term in Eq. \eqref{eq:zinonautgenef_NA} depends on the nonautonomous part of the piezoelectric stiffness. We notice that this last term is usually not diagonalised by the system eigenfunctions. Also, following~\cite{opreni22high}, it can be observed that the structure of both problems (autonomous and non-autonomous) are the same, such that a single problem can be rewritten, containing the two subcases.  Introducing the general superscript $\mathring{(\cdot)}$ to identify mappings and reduced dynamics of either autonomous and non-autonomous equations, one can rewrite the two previous problems as a single one which reads:

  \begin{align}\label{eq:first_homo_matrix_veps1_NA}
    & \forall\, \,\Is \in \cH^{(p)},\quad \forall\, p \in \{1,...,q\}, \nonumber \\
    & \left(\gsigma_{\Is}
    \left[\begin{array}{cc}
\Mu & \zerou \\
\zerou & \Mu
\end{array}\right]+
\left[\begin{array}{cc}
\Cu & \Ku_{T} \\
-\Mu & \zerou
\end{array}\right]\right)
\left[\begin{array}{l}
\gWV_{\Is}^{(p)} \\
\gWU_{\Is}^{(p)}
\end{array}\right] + 
\sum_{s=1}^{2 n} \gf_{s \Is}^{(p)}
\left[\begin{array}{cc}
\Mu & \zerou \\
\zerou & \Mu
\end{array}\right] \Yu_{s} =  \left[\begin{array}{c}
 \gE_{\Is}^{(p)} \\
-\Mu \gFU_{\Is}^{(p)}
\end{array}\right].
  \end{align}

where the excitation vector $\gE_{\Is}^{(p)}$ follows the following definition for the $\veps^{0}$ homological equations:

\begin{equation}\label{eq:defXiSigaut_NA}
 \gE_{\Is}^{(p)}   \delequal \Eu_{\Is}^{(p)} \delequal -\Mu\FV_{\Is}^{(p)} -{\FG}_{\Is}^{(p)}-{\FH}_{\Is}^{(p)}.
\end{equation}

On the other hand, the $\veps^{1}$ homological equations have an exitation vector defined as:

\begin{equation}\label{eq:defXiSigNONaut_NA}
  \gE^{(p)}_{\Is} \delequal \hat{\Eu}_{\Is}^{(p)} \delequal  -\Mu\hFV_{\Is}^{(p)} -\hFG_{\Is}^{(p)}-\hFH_{\Is}^{(p)} + \hFu^{(p)}_{\Is}.
\end{equation}

\section{Cantilever beam}
\label{app:cantilever}

In this Section we collect some results concerning the simulation of a cantilever beam
in large bending.
Although less meaningful from a technical point of view, 
this test is known to be a more severe benchmark than the clamped-clamped beam \cite{vizza21high,opreni22high} as it heavily involves non-resonant coupling between master and slave coordinates and can be employed to validate the limits of the proposed formulation.

\begin{figure}[p]
    \centering
    \includegraphics[width = .8\linewidth]{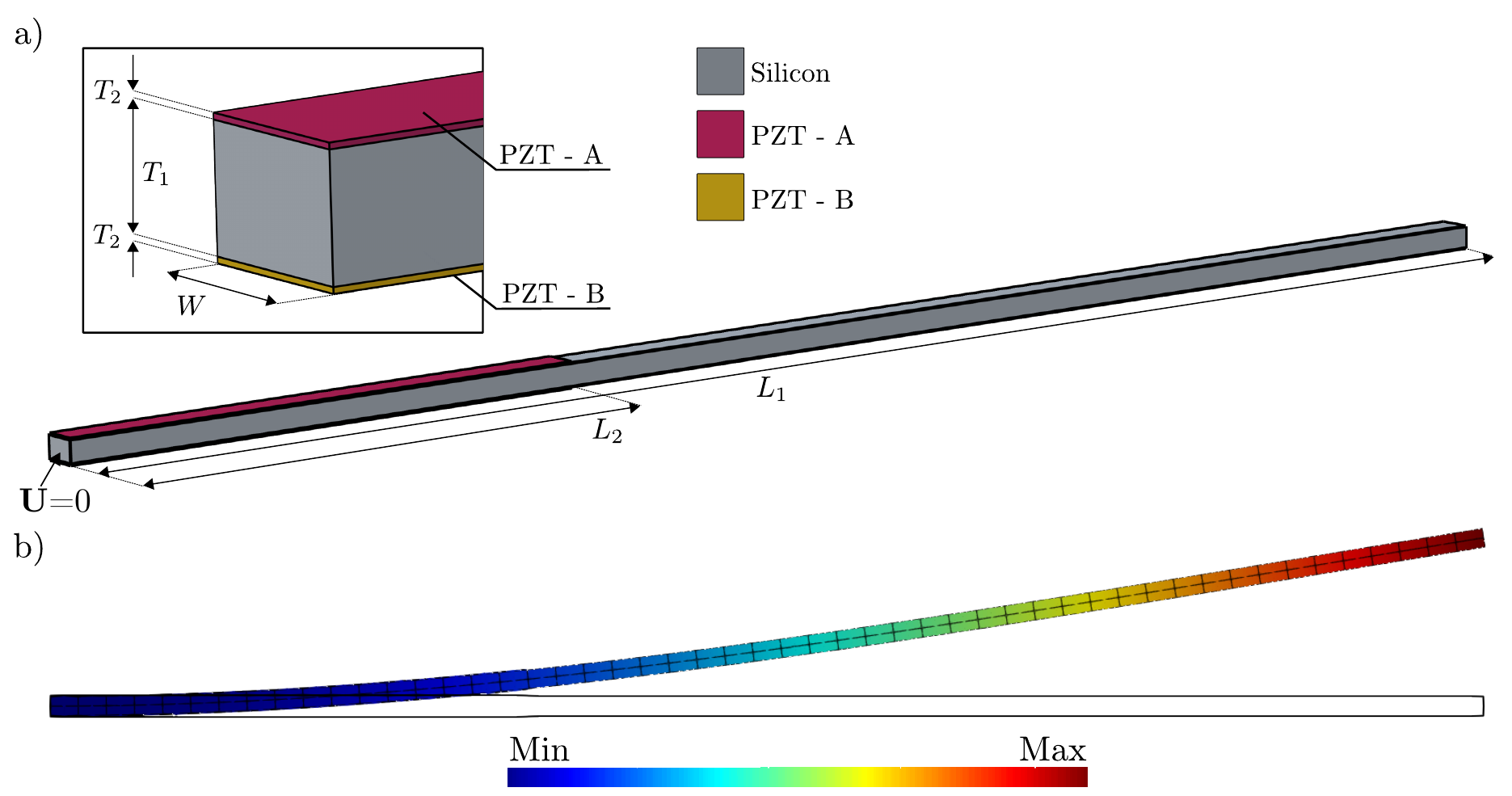}
    \caption{(a) Geometry of the tested cantilever. $T_1 = 0.2\,\mu$m, $T_2 = 0.01\,\mu$\,m, $L_1 = 15\,\mu$m, $L_2 = 5\,\mu$m, 
    $W = 0.3\,\mu$m. (b) Shape of the first bending mode $\bPhi_B$. Homogeneous Dirichlet boundary conditions are imposed on the initial cross section}
    \label{fig:cantilever_geo}
\end{figure}

The beam is illustrated in
Figure~\ref{fig:cantilever_geo}a, where $L_1 = 15\,\mu$m is the total length of the beam,
$L_2 = 5\,\mu$m is the length of the two piezo patches.
$T_1 = 0.2\,\mu$m denotes the thickness of the silicon body of the beam, while 
$T_2 = 0.01\,\mu$m is the piezo thickness.
The materials properties are detailed in Table \ref{tab:parameters_cc}.
As for the previous example, piezo patches are symmetrically positioned on the upper and lower surfaces of the beam,
but in this configuration, they cover one third of the beam length
at the anchoring point.
The two patches are actuated according to Eq.~\eqref{eq:actuation_law}
and set in resonant motion the first bending mode illustrated 
in Figure~\ref{fig:cantilever_geo}.

\begin{figure}[p]
    \centering
    \includegraphics[width = .99\linewidth]{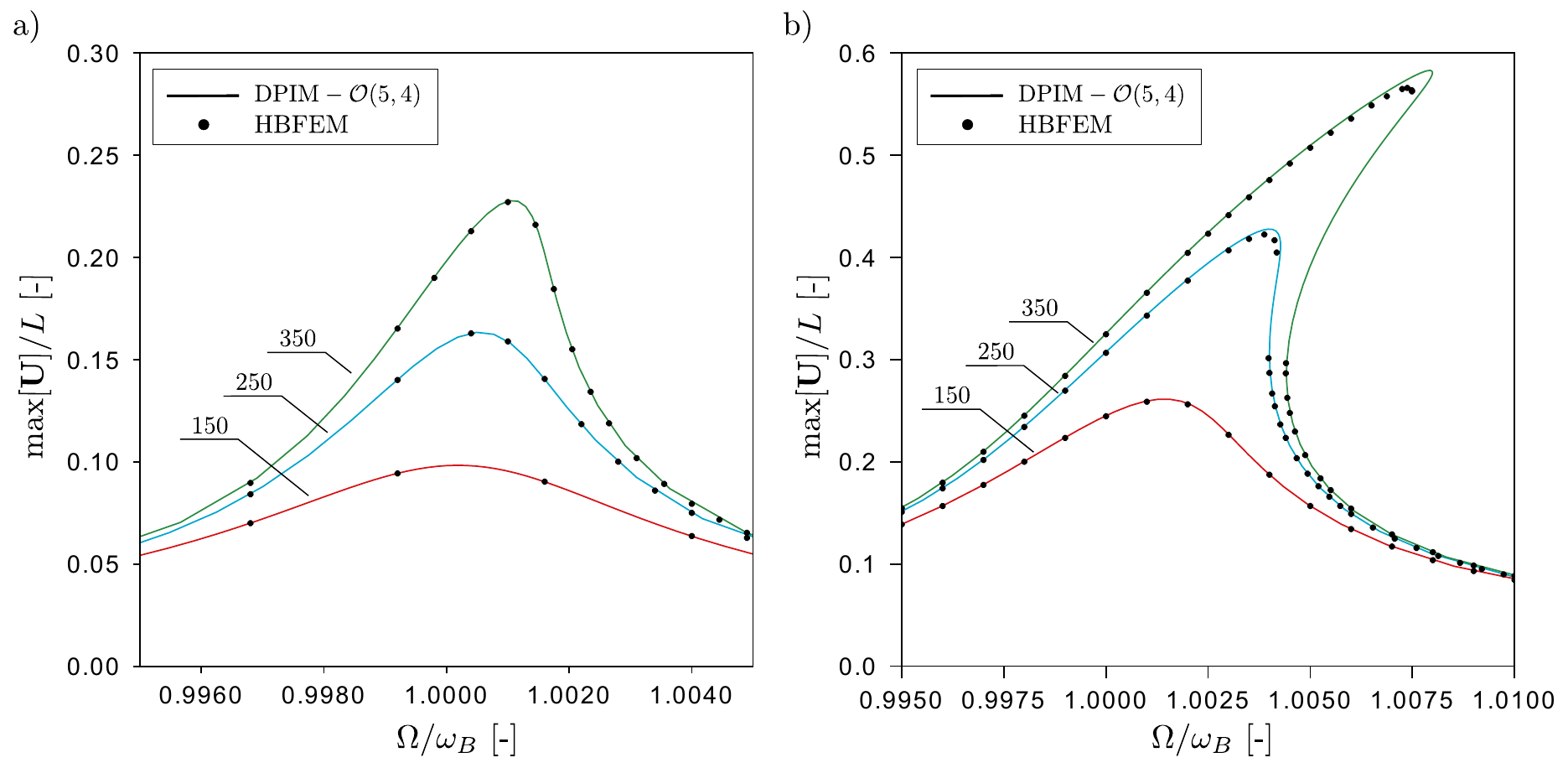}
    \caption{Comparison between full order HBFEM simulations and DPIM reduced model computed for 10 V (a) and 15 V (b). The tags in the charts report the quality factors $Q$ used for the analyses.}
    \label{fig:cantilever_frf}
\end{figure}

\begin{figure}[p]
    \centering
    \includegraphics[width = 0.7\linewidth]{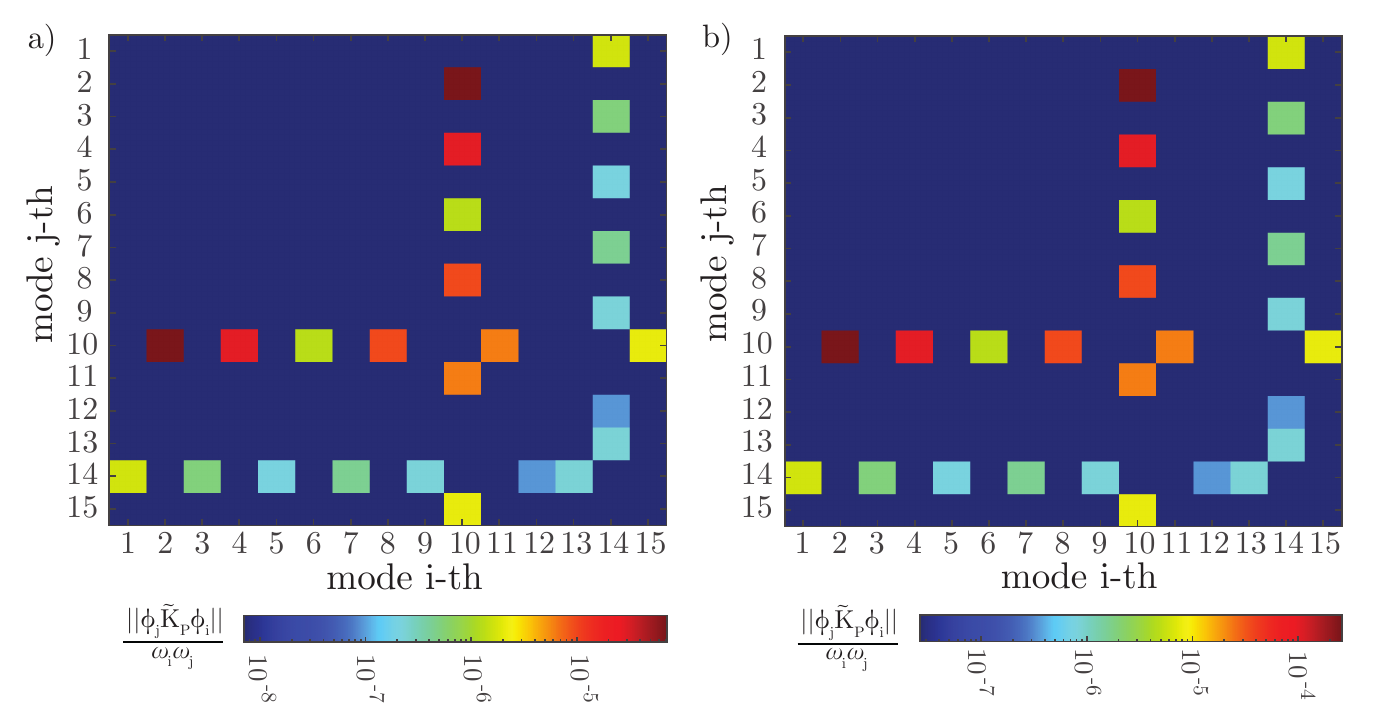}   
    \caption{Visual representation of the influence of $\tKu_P$ on the first 15 eigenmodes. 
    Fig. a) and b) refer to the voltage levels 10V, 15V corresponding to the polarisation curves of Fig.~\ref{fig:polarisation_curves}. Each coloured dot
    in the matrices
    represents  $|\bPhi^T_j\tKu_P\bPhi_i|/(\omega_i \omega_j)$}.
    \label{fig:kpiezo_cant}
\end{figure}

The polarisation curves reported in Fig.~\ref{fig:polarisation_curves}a)
have been used and a Rayleigh damping with $\beta=0$ and 
$\alpha$ calibrated so as to obtain the quality factors indicated in 
Fig.~\ref{fig:cantilever_frf}.
Figures~\ref{fig:cantilever_frf}a and \ref{fig:cantilever_frf}b collect the simulations for the two voltage bias
$V_0=10$\,V and $V_0=15$\,V, respectively. In each chart the analyses have been repeated for three different values of the quality factor $Q$, 
namely 150, 250 and 350.
It is worth stressing that the vertical axis reports the maximum displacement normalised with respect to the length of the beam,
so that the deviation evidenced between the DPIM results and the reference HBFEM 
occurs only at very large vibration amplitudes and may denote
the validity limit of the first-order expansion in $\varepsilon$ 
in Eqs.\eqref{eq:zeNLmaps},\eqref{eq:reduceddyn00} is being progressively 
reached.

As a final check, we verify that the linear effect of the matrix $\tKu_P$ is  still negligible, to ensure that this effect has not become important as compared to the clamped-clamped beam, which might give another explanation for the small discrepancies observed at very large amplitudes. To that purpose, the same analysis as at the end of Section~\ref{sec:CCbeam} is reproduced, by quantifying the magnitude of the terms brought by the matrix $\tKu_P$ in the modal space of the tangent problem corresponding to the new static position. The results are reported in Fig.~\ref{fig:kpiezo_cant}. The coupling pattern is different from the clamped beam case, and one can observe that the magnitudes are larger, with a maximum close to 3.10$^{-4}$ for the maximum voltage considered at 15~V. The coupling contribution is thus more meaningful, but still appears as small and should not be considered here as the main obstacle in the application of the method.

\end{document}